\documentclass[12pt,a4paper,reqno]{amsart}
\usepackage[english]{babel}
\usepackage{amssymb,latexsym,amsfonts,amsthm,upref,amsmath}
\usepackage{tikz}
\usepackage[margin=1in]{geometry}
\usepackage[colorlinks=false]{hyperref}
\usepackage{tikz}
\hypersetup{urlcolor=blue, linkcolor=red, citecolor=green, filecolor=magenta}
\usepackage{comment} 
 
\newtheoremstyle{prim}{}{}{\normalfont}{}{\bfseries}{}{ }{}
\theoremstyle{prim}
\newtheorem{ex}{Example}
\newtheoremstyle{stil}{}{}{\slshape}{}{\bfseries}{}{ }{}
\theoremstyle{stil}
\newtheorem{thm}{Theorem}[section]
\newtheoremstyle{defi}{}{}{}{}{\bfseries}{}{ }{}
\theoremstyle{defi}

\theoremstyle{defi}
\newtheorem{rem}[thm]{Remark}
\theoremstyle{stil}
\newtheorem{pro}[thm]{Proposition}
\theoremstyle{stil}
\newtheorem{lem}[thm]{Lemma}
\theoremstyle{stil}
\newtheorem{cor}[thm]{Corollary}
\newenvironment{prf}{\noindent \textit{Proof.}}{\null\hfill$\qed$\hskip 2mm\vskip 2mm}

\newcommand{\rez}{\mathop{\mathrm{Res}}}

\newcommand{\tr}{\mathop{\mathrm{tr}}}
\newcommand{\spn}{\mathop{\mathrm{span}}}
\newcommand{\ch}{\mathop{\mathrm{ch}}}
\newcommand{\lt}{\mathop{\mathrm{lt}}}

\numberwithin{equation}{section}

\begin{document}
\frenchspacing

\title[Quasi-particles in the principal picture of $\widehat{\mathfrak{sl}}_{2}$ and Rogers-Ramanujan-type identities]{Quasi-particles in the principal picture of $\widehat{\mathfrak{sl}}_{2}$ and Rogers-Ramanujan-type identities}

\author{Slaven Ko\v{z}i\'{c}} 
\address[Slaven Ko\v{z}i\'{c}]{Department of Mathematics, University of Zagreb, Zagreb, Croatia}
\email{kslaven@math.hr}

\author{Mirko Primc}
\address[Mirko Primc]{Department of Mathematics, University of Zagreb, Zagreb, Croatia}
\email{primc@math.hr}



\thanks{The research was supported by the grant ``Algebraic and combinatorial methods in vertex algebra theory'' of the Croatian Science Foundation.}

\begin{abstract} 
In their seminal work J. Lepowsky and R. L. Wilson gave a vertex-operator theoretic interpretation of Gordon-Andrews-Bressoud's generalization of 
Rogers-Rama\-nujan 
combinatorial identities, by constructing bases 
of vacuum spaces for the principal Heisenberg subalgebra
of standard $\widehat{\mathfrak{sl}}_{2}$-modules, parametrized with partitions satisfying
certain difference 2 conditions.
In this paper we define quasi-particles in the principal picture of  $\widehat{\mathfrak{sl}}_{2}$ and construct quasi-particle monomial bases
of standard $\widehat{\mathfrak{sl}}_{2}$-modules for which principally specialized characters are given as products of sum sides
of the corresponding analytic Rogers-Ramanujan-type identities with the character of the Fock space for the principal Heisenberg subalgebra.
\end{abstract}

\maketitle  
 

\tableofcontents

\section{Introduction} 
Famous Rogers-Ramanujan identities are two analytic identities, for $a=0,1$,
\begin{align*}
\prod_{m \geq 0}\frac{1}{(1-q^{5m+1+a})(1-q^{5m+4-a})}&= \sum_{m\geq0}\frac{q^{m^2+am}}{(1-q)(1-q^{2})\cdots (1-q^m)}.
\end{align*}
If, for $a=0$, we expand both sides in Taylor series, then the  coefficient of $q^n$ obtained from the product side can be interpreted 
as a number of partitions of $n$ with parts congruent $\pm 1\textrm{ mod }5$. On the other side, the $n$th coefficient of the summand
\begin{equation*}
\frac{q^{m^2}}{(q)_m},\qquad (q)_m =(1-q)(1-q^{2})\cdots (1-q^m),
\end{equation*}
can be interpreted as a number of partitions of $n$ into $m$ parts such that a difference between two consecutive parts is at least two.
Namely, $q^{m^2}$ ``represents'' the partition $(2m-1)+\ldots+3+1=m^2$, on which we ``add'' partitions 
$\lambda_1 +\ldots+\lambda_s$ with parts at most $m$  ``counted'' by $1/(q)_m$. For example, $q^{4^2}$ ``represents''
the ``smallest'' partition  $7+5+3+1=16$ satisfying difference condition, on which we ``add'' partitions
$\lambda_1 +\ldots+\lambda_s$ with parts at most $4$ ``counted'' by $1/(q)_4$. We can illustrate this construction for partition 
$12+9+7+2=30$ graphically,
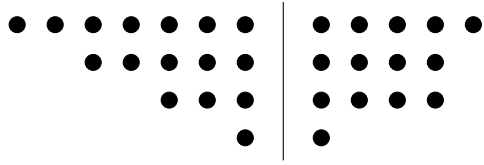
\begin{center}
\begin{tikzpicture}[scale=1]
\draw (0.5,-0.3) -- (0.5,1.8);
\filldraw (-3,1.5) circle (3pt);\filldraw (-2.5,1.5) circle (3pt);\filldraw (-2,1.5) circle (3pt);\filldraw (-1.5,1.5) circle (3pt);\filldraw (-1,1.5) circle (3pt);\filldraw (-0.5,1.5) circle (3pt);\filldraw (0,1.5) circle (3pt);\filldraw (1,1.5) circle (3pt);\filldraw (1.5,1.5) circle (3pt);\filldraw (2,1.5) circle (3pt);\filldraw (2.5,1.5) circle (3pt);\filldraw (3,1.5) circle (3pt);
\filldraw (-2,1) circle (3pt);\filldraw (-1.5,1) circle (3pt);\filldraw (-1,1) circle (3pt);\filldraw (-0.5,1) circle (3pt);\filldraw (0,1) circle (3pt);\filldraw (1,1) circle (3pt);\filldraw (1.5,1) circle (3pt);\filldraw (2,1) circle (3pt);\filldraw (2.5,1) circle (3pt);
\filldraw (-1,0.5) circle (3pt);\filldraw (-0.5,0.5) circle (3pt);\filldraw (0,0.5) circle (3pt);\filldraw (1,0.5) circle (3pt);\filldraw (1.5,0.5) circle (3pt);\filldraw (2,0.5) circle (3pt);\filldraw (2.5,0.5) circle (3pt);
\filldraw (0,0) circle (3pt);\filldraw (1,0) circle (3pt);
\end{tikzpicture}
\end{center}
 where we ``added'' partition $4+3+3+3+1=14$.

Analytic Rogers-Ramanujan identities have Gordon-Andrews-Bressoud's generalization
\begin{equation}\label{grr}
\prod_{\substack{n\geq 1\\n\not\equiv 0,\pm r(\textrm{mod }2l+s)}}(1-q^n)^{-1} = 
\sum_{n_1 , n_2 , \ldots ,  n_{l-1}\geq 0} \frac{q^{N_{1}^{2}+N_{2}^{2}+...+ N^{2}_{l-1}+N_{r}+N_{r+1}+\ldots+N_{l-1}   }}{(q)_{n_1 } (q)_{n_2 }\cdots (q)_{n_{l-2}} (q^{2-s})_{n_{l-1}}},
\end{equation}
where $N_j=n_j+n_{j+1}+\ldots+n_{l-1}$ and $s=0,1$, $l\geq 2$, $1\leq r\leq l-1$ (cf. \cite{Go}, \cite{A1}, \cite{A2}, \cite{Br1}, \cite{Br2}).
These identities have also a combinatorial interpretation, but for $2l+s>5$ it is not so easy to interpret the sum side of
(\ref{grr}) as a generating function for a number of partitions satisfying certain difference $2$ conditions among parts.

In 1980's Rogers-Ramanujan-type identities appeared in statistical physics and in representation theory of affine Kac-Moody Lie algebras.
This led to two lines of intensive research and numerous generalizations of both analytic and combinatorial identities,
the reader may consult, for example, the papers \cite{BM}, \cite{CLM}, \cite{JMS}, \cite{W} and the references therein.

J. Lepowsky and S. Milne discovered in \cite{LM} that the product sides of (\ref{grr}) multiplied with 
factor $F=\prod_{n\geq 1}(1-q^{2n-1})^{-1}$ are  the so-called principally specialized characters of standard modules for affine 
Lie algebra $\widehat{\mathfrak{sl}}_2$. 
Lepowsky and R. L. Wilson 
realized that the factor $F$ is a character of the Fock space for the principal Heisenberg subalgebra of  $\widehat{\mathfrak{sl}}_2$, 
and that the sum sides of (\ref{grr}) are the principally specialized characters of the vacuum spaces of standard modules for the action of principal Heisenberg subalgebra. 
In a series of papers (cf. \cite{LW3}--\cite{LW4}) Lepowsky and Wilson discovered vertex operators \emph{in the principal picture} 
on standard $\widehat{\mathfrak{sl}}_2$-modules and constructed bases 
of vacuum spaces for the principal Heisenberg subalgebra 
parametrized by partitions satisfying certain difference 2 conditions. 
Nowadays, we may say that combinatorial Rogers-Ramanujan identities come from two theories for $\widehat{\mathfrak{sl}}_2$: 
the product side from the representation theory for Kac-Moody Lie algebras and the sum side from the representation theory for vertex operator algebras.

Inspired by Lepowsky-Wilson's approach, in \cite{LP2} were constructed combinatorial bases for standard $\widehat{\mathfrak{sl}}_2$-modules by using
Frenkel-Kac-Segal vertex operators \emph{in the homogeneous picture}.
By using  results of Andrews, character formulas similar to the sum side of (\ref{grr}) were obtained in \cite{LP2}.
B. L. Feigin and A. V. Stoyanovsky introduced in \cite{FS} so-called quasi-particles, coefficients of formal Laurent series
\begin{equation*}
X_{\alpha}(z)^{p}=X_{\alpha}(z_1)X_{\alpha}(z_2)\cdots X_{\alpha}(z_p) \big|_{z_1=z_2=\ldots=z_p=z},
\end{equation*} 
in order to interpret
``directly'' the character formulas in \cite{LP2}.
Later on, G. Georgiev in \cite{Ge} constructed another type of combinatorial bases for standard   $\widehat{\mathfrak{sl}}_n$-modules---essentially 
the monomials in quasi-particles in the homogeneous picture---from which one can easily obtain  character formulas similar to the sum side of (\ref{grr})
by using the argument we have illustrated above.

In this paper we define quasi-particles in the principal picture of  $\widehat{\mathfrak{sl}}_{2}$ and construct quasi-particle monomial bases
of standard $\widehat{\mathfrak{sl}}_{2}$-modules for which principally specialized characters are given as products of sum sides of (\ref{grr})
with the  factor $F$. Similar approach was proposed by  Wilson in \cite{Wi} where the sum side of (\ref{grr}) is combinatorially related
to quasi-particle monomials in the principal picture.

Although  in the principal picture  the interesting identities ``live'' on the vacuum space, we work, as in  \cite{MP}, with the whole standard module, 
including the action of the principal Heisenberg subalgebra. 
In  \cite{MP}, the linear independence of Lepowsky-Wilson's combinatorial bases of standard $\widehat{\mathfrak{sl}}_2$-modules 
is proved by constructing the ``complementary'' bases in maximal submodules $W(\Lambda)$ of Verma $\widehat{\mathfrak{sl}}_2$-modules. 
Here we use the setting and the results of \cite{MP} to construct Georgiev-type quasi-particle monomial bases of 
standard $\widehat{\mathfrak{sl}}_2$-modules  in the principal picture.
As a result, it easily follows that the sum side of (\ref{grr}) is a 
natural 
factor of the character of a standard module, but only for some $\Lambda$ we are able to construct a ``complementary'' basis of $W(\Lambda)$ and, as a consequence, prove the linear independence of quasi-particle monomial basis of $L(\Lambda)$.

We start with  basis elements $X(n)$, $B(2n+1)$, $n\in\mathbb{Z}$ in the principal picture of $\widehat{\mathfrak{sl}}_2$ and
use formal Laurent series
\begin{equation}\label{i_rough}
X^{(p)}(\zeta)=\prod_{1\leq i<j\leq p}\left(\zeta_i\zeta_j\right)^{-1}\left(\zeta_i+\zeta_j\right)^{2}X(\zeta_1)\cdots X(\zeta_p)\Big|_{\zeta_1 =\zeta_2=\ldots=\zeta_p=\zeta}
\end{equation}
introduced in \cite{MP}. Since $X^{(p)}(\zeta)$ resembles the $p$-th power $X_{\alpha}(z)^p$ in the homogeneous picture, we call the coefficients
$$X^{(p)}(n),\quad n\in\mathbb{Z},\,p\geq 1,$$
of formal Laurent series $X^{(p)}(\zeta)$ {\em the quasi-particles of $\widehat{\mathfrak{sl}}_2$ in the principal picture}.
The obvious relations like
\begin{equation}\label{i_relacije}
 X_{\alpha}(z)^s X_{\alpha}(z)^p = X_{\alpha}(z)^{s-1} X_{\alpha}(z)^{p+1}
 \end{equation}
for quasi-particles in the homogeneous picture, have more complicated generalizations for quasi-particles in the principal picture (cf. Lemma \ref{lemma2}).
In Section \ref{5} we construct Georgiev-type quasi-particle bases of Verma $\widehat{\mathfrak{sl}}_2$-modules $M(\Lambda)$ 
(with a highest weight vector $v_{\Lambda}$),
\begin{equation}\label{i_baza}
B(i_1)\cdots B(i_r)X^{(p_1)}(j_1)\cdots X^{(p_s)}(j_s)v_{\Lambda},
\end{equation}
$r,s\geq 0$, for odd $i_1\leq \ldots\leq i_r\leq -1$ and for $1\leq p_1\leq \ldots\leq p_s$, 
satisfying ``difference conditions''
\begin{equation}\label{i_uvjet1}
j_l\leq -2p_l + j_{l+1} \quad\textrm{if}\quad p_l=p_{l+1}
\end{equation}
and ``initial conditions''
\begin{equation}\label{i_uvjet2}
j_l\leq -p_l -2p_{l}(s-l) \quad\textrm{if}\quad p_l<p_{l+1}
\end{equation}
(see Theorem \ref{verma:main}).
In fact, the above mentioned generalization of relations (\ref{i_relacije}) enables us to reduce the PBW spanning set of $M(\Lambda)$ to a
spanning set of the form (\ref{i_baza})--(\ref{i_uvjet2}), and the linear independence follows from Georgiev's construction of quasi-particle
bases of principal subspaces of standard $\widehat{\mathfrak{sl}}_2$-modules in \cite{Ge} (see  \cite{Bu}).

On a level $k$ standard module $L(\Lambda)$ we have relations
$$R_{p}(\zeta)=a_p X^{(p)}(\zeta)-(-1)^{k_0}a_q E^{-}(-\zeta)X^{(q)}(-\zeta)E^{+}(-\zeta)=0$$
for $p+q=k$, and as a consequence,
on the standard module $L(k\Lambda_0)=M(k\Lambda_0)/W(k\Lambda_0)$, $k$ odd, we have a Georgiev-type monomial bases of the form
(\ref{i_baza}) with the restriction
$$1\leq p_1\leq \ldots\leq p_s\leq [k/2]$$
(cf. Theorem \ref{standard:main}). Roughly speaking, the term
$$\frac{q^{p\,n_{p}^{2}}}{(q)_{n_{p}}}$$
in the right-hand side of character formula (\ref{grr}) appears, as in the Rogers-Ramanujan case, from a monomial factor
$$X^{(p)}(j_{s-n_{p}+1})\cdots X^{(p)}(j_s)$$
in (\ref{i_baza}) for which difference $2p$ conditions (\ref{i_uvjet1}) hold (see Section \ref{7}).

Finally, for some $\Lambda$ we construct ``complementary'' bases of $W(\Lambda)$  by using coefficients of $ R_{p}(\zeta)$ and quasi-particle monomials (see Theorem \ref{last:thm}) and, as a consequence, 
we have a vertex-operator theoretic proof of some of Andrews' analytic  identities (\ref{grr}).

Most of relations for $ R_{p}(\zeta)$ which we use are already proved in \cite{MP}, so we use all the notation from \cite{MP}.
After \cite{MP} was published, the theory of twisted modules for vertex operator algebras was developed (cf. \cite{L}), and we show that 
$X^{(p)}(\zeta)$ corresponds to iterate product 
\begin{equation}\label{i_voa}
\underbrace{x(z)_{-1}\ldots x(z)_{-1}}_{p-1}x(z)
\end{equation}
of the twisted field $x(z)$, which corresponds to $X(\zeta)$ with $\zeta=z^{-1/2}$ (see Remark \ref{onlyremark}).
However, we make no use of (\ref{i_voa}) since then there is
no need for ``translating'' the results in \cite{MP}, while
the ``rough'' construction (\ref{i_rough}) seems to be easy enough to handle.
And last, but not the least,
we did not extend the twisted vertex operator algebra theory setting needed to include the formal Laurent series  $ R_{p}(\zeta)$
which are used in \cite{MP} and in our construction, nor we replaced $R_{p}(\zeta)$ by some ``non-integrated'' elements which are 
expressed in terms of principal Heisenberg Lie algebra vertex operator $B(\zeta)$ instead of ``group'' elements $E^{\pm}(\zeta)$.

We should say that, in part, this work was motivated by a construction of quasi-particles for quantum affine 
algebra $U_q (\widehat{\mathfrak{sl}}_n)$ in \cite{Ko}.



\section{The principal picture of \texorpdfstring{$\widehat{\mathfrak{sl}}_{2}$}{sl2} and highest weight modules}\label{Sec:2}

In this section we recall the notation and results from \cite{MP}.
Let $\mathfrak{g}=\mathfrak{sl}_2$ with the standard basis
$$e=
\left(\begin{array}{rr}
0&1\\0&0
\end{array}\right),\qquad
f=
\left(\begin{array}{rr}
0&0\\1&0
\end{array}\right),\qquad
h=
\left(\begin{array}{rr}
1&0\\0&-1
\end{array}\right), 
$$
and invariant symmetric bilinear form
$$\left<x,y\right>=\tr xy.$$
Denote by $\mathbb{C}[t,t^{-1}]$ the algebra of Laurent polynomials in the indeterminate $t$.
Consider the involution of $\mathfrak{g}$ which is $1$ on $h$ and $-1$ on $e$ and $f$. 
The corresponding \emph{twisted affine Lie algebra} $\hat{\mathfrak{g}}=\widehat{\mathfrak{sl}}_{2}$ is the Lie algebra
$$\hat{\mathfrak{g}}=h\otimes\mathbb{C}[t^{2},t^{-2}]\oplus\spn\left\{e,f\right\}\otimes t\mathbb{C}[t^{2},t^{-2}]\oplus\mathbb{C}c\oplus\mathbb{C}d,$$
where $c$ is a nonzero central element in $\hat{\mathfrak{g}}$, and for $x,y\in\mathfrak{g}$ and $m,n\in\mathbb{Z}$,
\begin{align}
[x\otimes t^m,y\otimes t^n]&=[x,y]\otimes t^{m+n}+\left<x,y\right>\tfrac{m}{2}\delta_{m+n,0}c,\nonumber\\
[d,x\otimes t^m]&=mx\otimes t^m,\nonumber
\end{align}
(whenever $x\otimes t^m$ and $y\otimes t^ n$ are in $\hat{\mathfrak{g}}$).

Define
\begin{align}
&e_0=f\otimes t,&&e_1=e\otimes t,\nonumber\\
&f_0=e\otimes t^{-1},&&f_1=f\otimes t^{-1},\nonumber\\
&\textstyle h_0=-h\otimes 1+\frac{1}{2}c,&&h_1=h\otimes 1+\textstyle\frac{1}{2}c.\nonumber
\end{align}
Then the elements $e_i$, $f_i$, $h_i$, $i=0,1$, form a system of canonical generators of $\hat{\mathfrak{g}}=[\hat{\mathfrak{g}},\hat{\mathfrak{g}}]$,
viewed as the Kac-Moody Lie algebra with the Cartan matrix
$$A=\left(\begin{array}{rr}
2&-2\\-2&2
\end{array}\right).$$
Note that
$$c=h_0+h_1.$$

The relations
\begin{align}
&\deg e_i=-\deg f_i=1,&&i=0,1,\nonumber\\
&\deg h_i=\deg d=0,&&i=0,1,\nonumber
\end{align}
define a $\mathbb{Z}$-grading of $\hat{\mathfrak{g}}$ and $d$ is just the degree operator.  
This is the \emph{principal gradation} of the Lie algebra $\hat{\mathfrak{g}}$.

Set
$$B(n)=(e+f)\otimes t^n\quad\textrm{for }n\in 2\mathbb{Z}+1,$$
and
\begin{align}X(n)=\left\{\begin{array}{l@{\,\ }l}
 (f-e)\otimes t^n &\textrm{ if }n\in 2\mathbb{Z}+1,\\
 h\otimes t^n&\textrm{ if }n\in 2\mathbb{Z}.
\end{array}\right.\nonumber
\end{align}
Then the set
$$\left\{B(m),X(n),c,d : m\in 2\mathbb{Z}+1, n\in \mathbb{Z}\right\}$$
is a basis of $\hat{\mathfrak{g}}$. 
It is clear that the elements $B(n)$ and $X(n)$ have degree $n$, that is
\begin{align}
&[d,B(n)]=nB(n),\label{def:rel:dB}\\
&[d,X(n)]=nX(n).\label{def:rel:dX}
\end{align}
Furthermore we have
\begin{align}
&[B(m),B(n)]=m\delta_{m+n,0}c\quad\textrm{for }m,n\in 2\mathbb{Z}+1;\label{def:rel:1}\\
&[B(m),X(n)]=2X(m+n)\quad\textrm{for }m\in 2\mathbb{Z}+1,n\in\mathbb{Z};\label{def:rel:2}\\
&[X(m),X(n)]=(-1)^{m+1}2B(m+n)+(-1)^{m}m\delta_{m+n,0}c\quad\textrm{for }m,n\in\mathbb{Z},\label{def:rel:3}
\end{align}
and where by definition $B(m+n)=0$ for $m+n\in 2\mathbb{Z}$.

Set
\begin{align*}
&\mathfrak{s}_0=\mathbb{C}c+\mathbb{C}d,\\
&\mathfrak{s}_{\pm}=\coprod_{n>0}\mathbb{C}B(\pm n),\\
&\mathfrak{s}=\mathfrak{s}_{-}\oplus\mathfrak{s}_{0}\oplus\mathfrak{s}_{+}.
\end{align*}
We call $\mathfrak{s}_{-}\oplus\mathbb{C}c\oplus\mathfrak{s}_{+}$ \emph{principal Heisenberg subalgebra} and $\mathfrak{s}$ the 
\emph{extended principal Heisenberg subalgebra} of $\hat{\mathfrak{g}}$.

As in \cite{MP}, we also use the so-called {\em $\mathfrak{s}$-filtration} of the universal enveloping algebra $U=U(\hat{\mathfrak{g}})$ 
of $\hat{\mathfrak{g}}$ (cf. \cite{LW2}): For $j\in\mathbb{Z}$ set
\begin{align}
&U_{(j)}=(0)\quad\textrm{if }j<0,\label{s:filtration}\\
&U_{(0)}=U(\mathfrak{s})\nonumber
\end{align}
and if $j>0$ denote by $U_{(j)}$ the linear span of all elements $x_1,\ldots,x_n\in U$, where
$x_1,\ldots,x_n\in\hat{\mathfrak{g}}$ and at most $j$ of the elements $x_r$ lie outside the subalgebra $\mathfrak{s}$. Obviously
\begin{align*}
&(0)=U_{(-1)}\subset U_{(0)}\subset U_{(1)}\subset \ldots \subset U,\\
&U=\bigcup_{j\geq 0}U_{(j)}.
\end{align*}
It is obvious from the commutation relations (\ref{def:rel:2}) and (\ref{def:rel:3}) that for any permutation $\sigma$
and $m_1,\ldots,m_n\in\mathbb{Z}$ we have
\begin{equation}\label{filtration:eq}
X(m_{\sigma(1)})\cdots X(m_{\sigma(n)})- X(m_1)\cdots X(m_n)\in U_{(n-1)}.
\end{equation}
For $j$ a positive integer, set
$$\Xi_{(j)}=\spn\left\{X(m_1)\cdots X(m_n) : n\leq j, m_1\leq m_2\leq \ldots\leq m_n\right\}.$$
Clearly, (\ref{filtration:eq}) and the commutation relations (\ref{def:rel:1}) and (\ref{def:rel:2}) imply
\begin{equation}\label{U_(j):eq}
U_{(j)}=U(\mathfrak{s}_{-})\Xi_{(j)}U(\mathfrak{s}_{0}+\mathfrak{s}_{+}).
\end{equation}
If $i<j$, then we say that elements in $U_{(i)}$ are ``shorter'' in $\mathfrak{s}$-filtration then the elements in $U_{(j)}$.

Set
\begin{align}
&\mathfrak{h}=\mathbb{C}h\oplus\mathbb{C}c\oplus\mathbb{C}d,\nonumber\\
&\mathfrak{n}_{\pm}=\spn\left\{x\in\hat{\mathfrak{g}} : \pm\deg x>0\right\}.\nonumber
\end{align}
Then we have a triangular decomposition
$$\hat{\mathfrak{g}}=\mathfrak{n}_{-}\oplus\mathfrak{h}\oplus\mathfrak{n}_{+}$$
and the corresponding notion of  highest weight $\hat{\mathfrak{g}}$-modules.
We shall restrict our attention to the highest weight modules with the highest weight $\Lambda\in\mathfrak{h}^{*}$
such that $\Lambda(d)=0$. Then every such $\hat{\mathfrak{g}}$-module $V$ is $\mathbb{Z}$-graded,
\begin{equation*}
V=\coprod_{n\leq 0}V_n,\qquad V_{0}=\mathbb{C}v_0,
\end{equation*}
where $V_{n}$ is a finite-dimensional eigenspace of $d$ with eigenvalue $n\in\mathbb{Z}$.

We define the \emph{principally specialized character} of $V$ as the formal power series
$$\textstyle\ch_{q}V=\sum_{n\geq 0}(\dim V_{-n})q^{n}.$$

As a graded vector space, the Verma module $M(\Lambda)$ (with the highest weight $\Lambda$)
is isomorphic to the universal enveloping algebra $U(\mathfrak{n}_{-})$. Hence
\begin{equation}\label{verma:char}
\textstyle\ch_q M(\Lambda)=F\cdot \displaystyle\prod_{n\geq 1}(1-q^n)^{-1},
\end{equation}
where
$$F=\prod_{n\geq 1}(1-q^{2n-1})^{-1}.$$

\section{Quasi-particles in the principal picture}

We shall use formal Laurent series as  in \cite[Section 5]{MP}. In particular, we use Laurent series
\begin{align*}
X(\zeta)&=\sum_{n\in\mathbb{Z}}X(n)\zeta^n,\\
B(\zeta)&=\sum_{n\in 2\mathbb{Z}+1}B(n)\zeta^n
\end{align*}
with coefficients $X(n)$ and $B(n)$ in the universal enveloping algebra $U(\hat{\mathfrak{g}})$ of the Lie algebra $\hat{\mathfrak{g}}$. 
Then we may write relations (\ref{def:rel:dB}) and (\ref{def:rel:dX}) as
\begin{align}
[d,B(\zeta)]&=DB(\zeta),\label{def:rel:zeta:dB}\\
[d,X(\zeta)]&=DX(\zeta),\label{def:rel:zeta:dX}
\end{align}
and relations (\ref{def:rel:1})--(\ref{def:rel:3}) as
\begin{align}
&[B(\zeta),B(\xi)]=c\sum_{n\in 2\mathbb{Z}+1}n(\zeta/\xi)^{n},\label{def:rel:zeta:1}\\
&[B(\zeta),X(\xi)]=2X(\xi)\sum_{n\in 2\mathbb{Z}+1}(\zeta/\xi)^{n},\label{def:rel:zeta:2}\\
&[X(\zeta),X(\xi)]=-2B(\xi)\delta(-\zeta/\xi)+c(D\delta)(-\zeta/\xi),\label{def:rel:zeta:3}
\end{align}
where $D$ is the linear operator $D_{\zeta}=\zeta(d/d\zeta)$ and
\begin{align}
&\delta(\zeta)=\sum_{n\in\mathbb{Z}}\zeta^n,\nonumber\\
&D\delta(\zeta)=\sum_{n\in\mathbb{Z}}n\zeta^n.\nonumber
\end{align}

The coefficients $X(n)$ of $X(\zeta)$ do not commute and the square $X(\zeta)^2$, i.e.
$$\lim_{\xi\to\zeta} X(\xi)X(\zeta)=\left(X(\xi)X(\zeta)\right)\Big|_{\xi=\zeta}$$
is not defined. However, we have the relation 
$$(\zeta\xi)^{-1}(\zeta+\xi)^{2}X(\zeta)X(\xi)=(\zeta\xi)^{-1}(\zeta+\xi)^{2}X(\xi)X(\zeta)$$
(cf. (\ref{def:rel:zeta:3})) and the limit
$$X^{(2)}(\zeta)=\lim_{\xi\to \zeta}(\zeta\xi)^{-1}(\zeta+\xi)^{2}X(\zeta)X(\xi)$$
exists. Moreover, for a positive integer $p$ the limit
\begin{equation}\label{Xp}
X^{(p)}(\zeta)=\lim_{\zeta_{i}\to \zeta}\prod_{1\leq i<j\leq p}\left(\zeta_i\zeta_j\right)^{-1}\left(\zeta_i+\zeta_j\right)^{2}X(\zeta_1)\cdots X(\zeta_p),
\end{equation} 
exists (cf. \cite[Corollary 5.8]{MP}). Series (\ref{Xp}) is ``almost'' the $p$-th power of the Laurent series $X(\zeta)$,
in a way parallel to $X_{\alpha}(\zeta)^{p}$ in the homogeneous picture (cf. \cite{FS}, \cite{Ge}, \cite{LP2}).
By following B. Feigin and A. Stoyanovsky we call $X^{(p)}(n)$, a coefficient in
$$X^{(p)}(\zeta)=\sum_{n\in\mathbb{Z}}X^{(p)}(n)\zeta^n,$$
the \emph{quasi-particle} of \emph{degree} $n$ and \emph{charge} $p$. With a proper interpretation  (cf. \cite{MP}), the quasi-particle
$X^{(p)}(n)$ is an infinite sum
\begin{equation*}
X^{(p)}(n)=\sum_{i_1+\ldots+i_p=n}f(i_1,\ldots,i_p)
\end{equation*}
of the summable family
$$\left\{f(i_1,\ldots,i_p) : i_1+\ldots+i_p=n\right\}$$
of coefficients in formal Laurent series
$$\prod_{1\leq i<j\leq p}\left(\zeta_i\zeta_j\right)^{-1}\left(\zeta_i+\zeta_j\right)^{2}X(\zeta_1)\cdots X(\zeta_p)=\sum_{i_1,\ldots,i_p\in\mathbb{Z}}f(i_1,\ldots,i_p)\zeta_{1}^{i_1}\cdots\zeta_{p}^{i_p}.$$
For an arbitrary positive integer $p$ the series $X^{(p)}(\zeta)$ satisfy recursive formula (cf. \cite[Lemma 5.12.]{MP})
\begin{equation}\label{Xp:recursive}
X^{(p+1)}(\zeta)=2^{2(p-1)}\lim_{\zeta_1,\zeta_2\to\zeta}(\zeta_1\zeta_2)^{-1}(\zeta_1+\zeta_2)^{2}X(\zeta_1)X^{(p)}(\zeta_2).
\end{equation}
Note that
$$\lim_{\zeta_{i}\to \zeta}\prod_{1\leq i<j\leq p}\left(\zeta_i\zeta_j\right)^{-1}\left(\zeta_i+\zeta_j\right)^{2}=2^{p(p-1)}.$$
The next proposition is an immediate consequence of (\ref{def:rel:zeta:2}).
\begin{pro}\label{pro:Xp:B}
For $p\in\mathbb{N}$ we have
$$[B(\zeta),X^{(p)}(\xi)]=2pX^{(p)}(\xi)\sum_{n\in 2\mathbb{Z}+1}(\zeta/\xi)^n .$$
\end{pro}

\begin{rem}\label{onlyremark}
In this remark we   derive series $X^{(p)}(\zeta)$, defined by (\ref{Xp}), by using products of twisted vertex operators 
introduced by H.-S. Li in \cite{L}.

Set
$$X_{0}(\zeta)=\sum_{n\in 2\mathbb{Z}}X(n)\zeta^n,\qquad X_{1}(\zeta)=\sum_{n\in 2\mathbb{Z}+1}X(n)\zeta^n.$$
Obviously $X(\zeta)=X_{0}(\zeta)+X_{1}(\zeta)$.
Define
\begin{align*}
&x_{0}(z)=\sum_{n\in 2\mathbb{Z}}x_{0}(n/2)z^{-\frac{n}{2}-1}=\sum_{n\in 2\mathbb{Z}}X(n)z^{-\frac{n}{2}-1},\\
&x_{1}(z)=\sum_{n\in 2\mathbb{Z}+1}x_{1}(n/2)z^{-\frac{n}{2}-1}=\sum_{n\in 2\mathbb{Z}+1}X(n)z^{-\frac{n}{2}-1},\\
&x(z)=\sum_{n\in \mathbb{Z}}x(n/2)z^{-\frac{n}{2}-1}=x_{0}(z)+x_{1}(z).
\end{align*}
Relation  (\ref{def:rel:zeta:3}) implies
\begin{align}
&(z_1-z)^{2}[x_{0}(z_1),x(z)]=(z_1-z)^{2}[x_{1}(z_1),x(z)]=0.\label{locality}
\end{align}

By applying the multiplication formula for twisted vertex operators, introduced in \cite{L}, we get
{\allowdisplaybreaks\begin{align*}
&Y(x_{0}(z),z_0)x(z)=\sum_{n\in\mathbb{Z}}(x_{0}(z)_{n}x(z))z_{0}^{-n-1}\\
=&\rez_{z_1}\left(z_{0}^{-1}\delta\left(\frac{z_1 -z}{z_0}\right)x_{0}(z_1)x(z)-z_{0}^{-1}\delta\left(\frac{z -z_1}{-z_0}\right)x(z)x_{0}(z_1)\right)\\
=&\rez_{z_1}\left(z_{0}^{-1}\delta\left(\frac{z_1 -z}{z_0}\right)x_{0}(z_1)x(z)-z_{0}^{-1}\delta\left(\frac{z -z_1}{-z_0}\right)\frac{(z - z_1)^{2}}{(-z_{0})^{2}}x(z)x_{0}(z_1)\right),\\[1em]
&Y(x_{1}(z),z_0)x(z)=\sum_{n\in\mathbb{Z}}(x_{1}(z)_{n}x(z))z_{0}^{-n-1}\\
=&\rez_{z_1}\left(\frac{z_1 - z_0}{z}\right)^{1/2}\left(z_{0}^{-1}\delta\left(\frac{z_1 -z}{z_0}\right)x_{1}(z_1)x(z)-z_{0}^{-1}\delta\left(\frac{z -z_1}{-z_0}\right)x(z)x_{1}(z_1)\right)\\
=&\rez_{z_1}\left(\frac{z_1 - z_0}{z}\right)^{1/2}\left(z_{0}^{-1}\delta\left(\frac{z_1 -z}{z_0}\right)x_{1}(z_1)x(z)-z_{0}^{-1}\delta\left(\frac{z -z_1}{-z_0}\right)\frac{(z - z_1)^{2}}{(-z_{0})^{2}}x(z)x_{1}(z_1)\right).
\end{align*}}
The expressions above can be simplified by using locality (\ref{locality}) and some properties of the classical 
delta function (cf. \cite[Proposition 2.3.8]{LL}):
\begin{align*}
&Y(x_{0}(z),z_0)x(z)=\rez_{z_1} z_{1}^{-1}\delta\left(\frac{z+z_0}{z_1}\right)x_{0}(z_1)x(z),\\
&Y(x_{1}(z),z_0)x(z)=\rez_{z_1}\left(\frac{z_1 - z_0}{z}\right)^{1/2} z_{1}^{-1}\delta\left(\frac{z+z_0}{z_1}\right)x_{1}(z_1)x(z).
\end{align*}

We will now restrict ourselves to calculating $x_{0}(z)_{-1}x(z)$ and $x_{1}(z)_{-1}x(z)$:
\begin{align}
x_{0}(z)_{-1}x(z)&=\rez_{z_0}\rez_{z_1} z_{0}^{-1}z_{1}^{-1}\delta\left(\frac{z+z_0}{z_1}\right)x_{0}(z_1)x(z)\nonumber\\
&=\rez_{z_0}\rez_{z_1} z_{0}^{-3}z_{1}^{-1}\delta\left(\frac{z+z_0}{z_1}\right)(z_1 -z)^{2}x_{0}(z_1)x(z)\nonumber\\
&=\rez_{z_0}\left( z_{0}^{-3} \left((z_1 -z)^{2}x_{0}(z_1)x(z)\right)\Big|_{z_1=z+z_0}\right).\label{rem:1}
\end{align}
Similarly, we get
\begin{equation}\label{rem:2}
x_{1}(z)_{-1}x(z)=\rez_{z_0}\left( z_{0}^{-3}\left((z_1 -z)^{2}(z_{1}/z)^{1/2}x_{1}(z_1)x(z)\right)\Big|_{z_1=z+z_0}\right).
\end{equation}
By introducing the series $X_{0}(\zeta)$ and $X(\zeta)$ in (\ref{rem:1}) we get
\begin{align*}
&x_{0}(z)_{-1}x(z)=\rez_{z_0}\left( z_{0}^{-3} \left((z_{1}z)^{-1}(z_1 -z)^{2}X_{0}(z_1^{-1/2})X(z^{-1/2})\right)\Big|_{z_1=z+z_0}\right)\\
=&\rez_{z_0}\left( z_{0}^{-3} \left((z_{1}z)^{-1}(z_1^{1/2} -z^{1/2})^{2}(z_1^{1/2} +z^{1/2})^{2}X_{0}(z_1^{-1/2})X(z^{-1/2})\right)\Big|_{z_1=z+z_0}\right)\\
=&\rez_{z_0}\left( z_{0}^{-3} \left(\underbrace{\left(\frac{(z_1^{1/2} -z^{1/2})^{2}}{(z_{1}z)^{1/2}}\right)}_{A(z_1,z)}\underbrace{\left(\frac{(z_1^{-1/2} +z^{-1/2})^{2}}{(z_{1}z)^{-1/2}}X_{0}(z_1^{-1/2})X(z^{-1/2})\right)}_{B(z_1,z)}\right)\Bigg|_{z_1=z+z_0}\right).
\end{align*}
Notice that
$$\left(A( z_1,z)B(z_1,z)\right)\Big|_{z_1=z+z_0}=A( z_1,z)\Big|_{z_1=z+z_0}B(z_1,z)\Big|_{z_1=z+z_0}.$$
By applying the binomial theorem expansion we get
$$A( z_1,z)\Big|_{z_1=z+z_0}=\frac{1}{4}\left(\frac{z_0}{z}\right)^{2}\left(1+p(z_0/z)\right)$$
for some polynomial $p(z)\in\mathbb{C}[z]$, $p(0)=0$. Since $B(z_1,z)\Big|_{z_1=z+z_0}$ contains only nonnegative powers of $z_0$, 
we finally conclude
\begin{equation}\label{rem:3}
x_{0}(z)_{-1}x(z)=\frac{1}{4z^{2}}\left(\frac{(z_1^{-1/2} +z^{-1/2})^{2}}{(z_{1}z)^{-1/2}}X_{0}(z_1^{-1/2})X(z^{-1/2})\right)\Bigg|_{z_1=z}.
\end{equation}
In a similar way we  obtain
\begin{equation}\label{rem:4}
x_{1}(z)_{-1}x(z)=\frac{1}{4z^{2}}\left(\frac{(z_1^{-1/2} +z^{-1/2})^{2}}{(z_{1}z)^{-1/2}}X_{1}(z_1^{-1/2})X(z^{-1/2})\right)\Bigg|_{z_1=z}.
\end{equation}
Since $$B( z_1,z)\Big|_{z_1=z}=X^{(2)}(z^{-1/2})$$ (recall (\ref{Xp})), equations (\ref{rem:3}) and (\ref{rem:4}) imply
\begin{equation*}
X^{(2)}(z^{-1/2})=4z^{2}x(z)_{-1}x(z),
\end{equation*}
where $x(z)_{-1}=x_{0}(z)_{-1}+x_{1}(z)_{-1}$.

The formula above  can be easily generalized to the higher charges. For an arbitrary positive integer $p$ we have
$$X^{(p)}(z^{-1/2})=2^{2(p-1)}z^{p}\underbrace{x(z)_{-1}\ldots x(z)_{-1}}_{p-1}x(z).$$ 

\end{rem}


\section{Relations for quasi-particles}

From now on it will be convenient to use definition (\ref{Xp}). 
Fix nonnegative integers $n_1,\ldots,n_p$ and let $D_{\zeta}=\zeta(d/d\zeta)$.
Then derivations of \ref{def:rel:zeta:3} and properties of $\delta$-function imply that limit
\begin{equation}
\lim_{\zeta_{i}\to\zeta}\prod_{1\leq i<j\leq p} (\zeta_i \zeta_j)^{-(N+1)}(\zeta_{i}+\zeta_{j})^{2(N+1)}
D_{\zeta_1}^{n_1}X(\zeta_1)\cdots D_{\zeta_p}^{n_p}X(\zeta_p)
\end{equation}
exists for $N=\max\left\{n_1,n_2,\ldots,n_p\right\}$.

\begin{lem}\label{lemma1}
Let $P(\zeta_1,\ldots,\zeta_p)$ be a Laurent polynomial such that the limit $\zeta_{i}\to\zeta$ exists for
$$C(\zeta_1,\ldots,\zeta_p):=P(\zeta_1,\ldots,\zeta_p)X(\zeta_1)\cdots X(\zeta_p).$$
Then for $r\geq 0$ we have
$$\lim_{\zeta_i\to\zeta} (D_{\zeta_1}+\ldots+D_{\zeta_p})^{r}C(\zeta_1,\ldots,\zeta_p)=D_{\zeta}^{r}\lim_{\zeta_i\to\zeta} C(\zeta_1,\ldots,\zeta_p).$$
\end{lem}

\begin{prf}
Set
$$C(\zeta_1,\ldots,\zeta_p)=\sum_{j_1,\ldots,j_p\in\mathbb{Z}}C(j_1,\ldots,j_p)\zeta_{1}^{j_1}\cdots\zeta_{p}^{j_p}.$$
By assumption, for every $n\in\mathbb{Z}$ the family of coefficients
$$\left\{C(j_1,\ldots,j_p) : j_1 +\ldots + j_p =n\right\}$$
is summable, i.e. for a given highest weight $\hat{\mathfrak{g}}$-module $V$ and a vector $v\in V$ there exists $N_{v,n}\in\mathbb{Z}$
such that
$$C(j_1,\ldots,j_p)v\neq 0\quad\textrm{implies}\quad j_1,\ldots,j_p\leq N_{v,n}$$
and hence
$$\sum_{j_1,\ldots,j_p\in\mathbb{Z}}C(j_1,\ldots,j_p)v\in V$$
is a finite sum of nonzero vectors. Hence we have
\begin{align*}
&\lim_{\zeta_i\to\zeta} (D_{\zeta_1}+\ldots+D_{\zeta_p})^{r}C(\zeta_1,\ldots,\zeta_p)\\
=&\lim_{\zeta_i\to\zeta} (D_{\zeta_1}+\ldots+D_{\zeta_p})^{r}\sum_{j_1,\ldots,j_p\in\mathbb{Z}}C(j_1,\ldots,j_p)\zeta_{1}^{j_1}\cdots\zeta_{p}^{j_p}\\
=&\lim_{\zeta_i\to\zeta} \sum_{j_1,\ldots,j_p\in\mathbb{Z}} (j_1 +\ldots  + j_p)^{r} C(j_1,\ldots,j_p)\zeta_{1}^{j_1}\cdots\zeta_{p}^{j_p}\\
=&\sum_{j_1,\ldots,j_p\in\mathbb{Z}} (j_1 +\ldots  + j_p)^{r}C(j_1,\ldots,j_p)\zeta^{j_1 + \ldots +j_p}\\
=&D_{\zeta}^{r}\sum_{j_1,\ldots,j_p\in\mathbb{Z}}C(j_1,\ldots,j_p)\zeta^{j_1 + \ldots +j_p}\\
=&D_{\zeta}^{r}\lim_{\zeta_i\to\zeta} C(\zeta_1,\ldots,\zeta_p).
\end{align*}
\end{prf}

\begin{lem}\label{lemma2}
Let $P(\zeta_1,\ldots,\zeta_s;\xi_1,\ldots,\xi_p)$ be a Laurent polynomial, $s\leq p$, such that there exists a limit
$$\lim_{\zeta_i,\xi_j\to\chi} P(\zeta_1,\ldots,\zeta_s;\xi_1,\ldots,\xi_p) \left(\prod_{i=1}^{s}D_{\zeta_i}^{n_i}X(\zeta_i)\right) \left(\prod_{j=1}^{p}X(\xi_j)\right)$$
for $n_1+\ldots+n_s\leq 2s-1$, and that for some Laurent polynomial $R(\zeta,\xi)$ we have
$$\lim_{\substack{\zeta_i\to\zeta\\\xi_j\to\xi}}P(\zeta_1,\ldots,\zeta_s;\xi_1,\ldots,\xi_p) \left(\prod_{i=1}^{s}X(\zeta_i)\right) \left(\prod_{j=1}^{p}X(\xi_j)\right)
=R(\zeta,\xi)X^{(s)}(\zeta)X^{(p)}(\xi).$$
Then for $r=0,1,\ldots,2s-1$
\begin{align}
\lim_{\zeta,\xi\to\chi}D_{\zeta}^{r}\left(R(\zeta,\xi)X^{(s)}(\zeta)X^{(p)}(\xi)\right)&\label{lemma2:1}\\
\equiv \left(X^{(p+1)}(\chi)\right),&\label{lemma2:2}
\end{align}
where (\ref{lemma2:2}) means that the coefficients of the formal Laurent series (\ref{lemma2:1}) in indeterminate $\chi$ are (``infinite'')
linear combinations of products of at most $s-1$ factors $X(n)$ with at least one coefficient of $X^{(p+1)}(\chi)$.
\end{lem}

\begin{prf}
We start with the formal Laurent series
\begin{equation}\label{lemma2:3}
(D_{\zeta_1}+\ldots+D_{\zeta_s})^{r} P(\zeta_1,\ldots,\zeta_s;\xi_1,\ldots,\xi_p) \prod_{i=1}^{s}X(\zeta_i) \prod_{j=1}^{p}X(\xi_j)
\end{equation}
 and evaluate the limit $\zeta_i,\xi_j\to\chi$ in stages in two ways:
 
 \noindent(1) First, by taking a limit $\zeta_i\to\zeta$ and using Lemma \ref{lemma1} we obtain
 \begin{equation*}
 D_{\zeta}^{r} \lim_{\xi_j\to\xi} R(\zeta;\xi_1,\ldots,\xi_p)X^{(s)}(\zeta)\prod_{j=1}^{p}X(\xi_j)
 \end{equation*}
 for some Laurent polynomial $R(\zeta;\xi_1,\ldots,\xi_p)$. Then, by taking the limit $\zeta,\xi_j\to\chi$ we obtain (\ref{lemma2:1}).
 
 \noindent(2) For $0\leq r\leq 2s-1$ we have
 \begin{equation}\label{s>p:1}
 (D_{\zeta_1}+\ldots+D_{\zeta_s})^{r}=\sum_{i_1,\ldots,i_r}D_{\zeta_{i_1}}\cdots D_{\zeta_{i_r}},
 \end{equation}
 where for each choice of $(i_1,\ldots,i_r)$ there is an index $t\in\left\{1,\ldots,s\right\}$ such that there is at most one $D_{\zeta_t}$  in the product
 \begin{equation}\label{s>p:2}
 D_{\zeta_{i_1}}\cdots D_{\zeta_{i_r}}.
 \end{equation}
 Then (\ref{lemma2:3}) can be written as a sum over $i_1,\ldots,i_r$. Any such summand can be written as
 $$\left(\prod_{i\neq t}D_{\zeta_i}^{a_i}\right)D_{\zeta_t}^{a}P(\zeta_1,\ldots,\zeta_s;\xi_1,\ldots,\xi_p) \prod_{i=1}^{s}X(\zeta_i) \prod_{j=1}^{p}X(\xi_j)$$
with $a\in\left\{0,1\right\}$. By taking a limit $\zeta_t,\xi_j\to\chi$ we obtain a factor
$D_{\chi}X^{(p+1)}(\chi)$ and/or $X^{(p+1)}(\chi)$ (cf. \cite[Lemma 5.12]{MP}). After that, by taking a limit $\zeta_i\to\chi$, $i\neq t$,  we obtain 
a series for which (\ref{lemma2:2}) holds.
\end{prf}

\begin{lem}\label{lemma3}
Let $N$ and $s$ be positive integers, $N>2s$. Let $D_{z}=z(d/d z)$ and $R(z)=(1+z)^{N}$.
Let $j$ be an integer. Set
$$a_{r,m}=D_{z}^{r}\left(z^{m}R(z)\right)\Big|_{z=1}.$$
Then $2s\times 2s$ matrix
\begin{equation*}
\left(a_{r,m}\right)_{r=0,1,\ldots,2s-1;\hspace{2pt}m=j,j+1,\ldots,j+2s-1},
\end{equation*}
is regular.
\end{lem}

\begin{prf}
For $r=0$ we have
$$\left(z^{m}(1+z)^{N}\right)\Big|_{z=1}=2^{N}$$
so the first row in our matrix is
$$(2^N,2^N,\ldots,2^N)=2^{N}(1,1,\ldots,1).$$
Write
$$(j,j+1,\ldots,j+2s-1)=(j_1,j_2,\ldots,j_{2s}).$$
For $r=1$ we have
$$D_{z}\left(z^{m}(1+z)^{N}\right)=mz^{m}(1+z)^{N}+Nz^{m+1}(1+z)^{N-1},$$
so the second row in our matrix is
\begin{align*}
&(j_1 2^N +N2^{N-1},j_2 2^N+N2^{N-1},\ldots, j_{2s} 2^N+N2^{N-1})\\
=& 2^N (j_1,j_2,\ldots,j_{2s}) + N2^{N-1}(1,1,\ldots,1).
\end{align*}
It is easy to see that
$$D_{z}^{r}\left(z^{m}(1+z)^{N}\right)=m^{r}z^{m}(1+z)^{N}+\sum_{l=1}^{r}f_{l}^{(r)}(m)z^{m+l}(1+z)^{N-l},$$
where $f_{l}^{(r)}(m)$ is a polynomial in $m$ of degree at most $r-l<r$. So, by using elementary transformations of rows,
our matrix $\left(a_{r,m}\right)$ can be reduced to the matrix
$$\left(\begin{array}{cccc}
1&1&\ldots&1\\
j_1&j_2&\ldots&j_{2s}\\
\vdots&\vdots&&\vdots\\
j_{1}^{2s-1}&j_{2}^{2s-1}&\ldots&j_{2s}^{2s-1}
\end{array}\right).
$$
\end{prf}

\begin{lem}\label{lemma4}
Let $s$ and $p$ be positive integers. Let $P(z)$ be a Laurent polynomial such that
$$P(\zeta/\xi)X^{(s)}(\zeta)X^{(p)}(\xi)=P(\zeta/\xi)X^{(p)}(\xi)X^{(s)}(\zeta).$$
Set
$$C(\chi)=\sum_{n\in\mathbb{Z}}C(n)\chi^{n}=\lim_{\zeta,\xi\to\chi}P(\zeta/\xi)X^{(s)}(\zeta)X^{(p)}(\xi).$$
Fix integers $n$ and $r_1,\ldots,r_q$, $q\geq 1$. Then for any given vector $v$ in any given highest weight $\hat{\mathfrak{g}}$-module $V$ 
there exists a finite subset $\left\{r_1,\ldots,r_q\right\}\subseteq J\subset \mathbb{Z}$ such that
a vector $C(n)v$ can be written as a finite sum
$$C(n)v=\sum_{r\in J}b_r X^{(s)}(r)X^{(p)}(n-r)v,$$
where $b_r=P(1)$ for $r=r_1,\ldots,r_q$.
\end{lem}

\begin{prf}
Let
$$P(z)=\sum_{l=-t}^{t}a_l z^l$$
and set
$$P(\zeta/\xi)X^{(s)}(\zeta)X^{(p)}(\xi)=\sum_{i,j\in\mathbb{Z}}f(i,j)\zeta^{i}\xi^{j}.$$
Then
\begin{align*}
f(i,j)&=\sum_{l=-t}^{t} a_l X^{(s)}(i-l)X^{(p)}(j+l)\\
&=\sum_{l=-t}^{t} a_l X^{(p)}(j+l)X^{(s)}(i-l).
\end{align*}
For a given vector $v$ in a given highest weight $\hat{\mathfrak{g}}$-module $V$ there exists $M>0$ such that
$$f(i,j)v=0\quad\textrm{for}\quad i\geq M \textrm{ or } j\geq M.$$
Hence the family
$$\left\{f(i,j)v : i+j=n\right\}$$
is summable and we can take a finite set $J$ such that
$$r_{1}+l,r_{2}+l,\ldots,r_{q}+l\in J\qquad\textrm{for all }l=-t,-t+1,\ldots,t$$
and such that
\begin{align*}
C(n)v&=\sum_{r\in J}f(r,n-r)
=\sum_{r\in J}\left(\sum_{l=-t}^{t}a_l X^{(s)}(r-l)X^{(p)}(n-r+l)v\right)\\
&=\sum_{l=-t}^{t}a_l \sum_{r'+l\hspace{1pt}\in\hspace{1pt} J} X^{(s)}(r')X^{(p)}(n-r')v.
\end{align*}
For $r=r_1,r_2,\ldots,r_q$ the coefficient $b_r$ of $X^{(s)}(r)X^{(p)}(n-r)v$ in the sum above equals
$$\sum_{t=-l}^{l}a_l=P(1).$$
\end{prf}

\begin{lem}\label{lemma5}
Let $s\leq p$ and $n,j\in\mathbb{Z}$. For any given vector $v$ in any given highest weight $\hat{\mathfrak{g}}$-module $V$ 
denote by $S^{(s,p)}_{j,n}v$ a sequence of $2s$ vectors
\begin{align}
&X^{(s)}(j)X^{(p)}(n-j)v,\nonumber\\
&X^{(s)}(j+1)X^{(p)}(n-(j+1))v,\label{lemma5:summands}\\
&\qquad\qquad\qquad\vdots\nonumber\\
&X^{(s)}(j+2s-1)X^{(p)}(n-(j+2s-1))v\nonumber
\end{align}
in a formal Laurent series $X^{(s)}(\zeta)X^{(p)}(\xi)v$. Then each vector in $S^{(s,p)}_{j,n}v$
can be written as a linear combination of vectors
$$X^{(s)}(i)X^{(p)}(n-i)v\notin S^{(s,p)}_{j,n}v$$
and vectors obtained by the action of monomials, that have  a coefficient of $X^{(p+1)}(\chi)$ as a factor, on $v$.
\end{lem}

\begin{prf}
We prove the statement of the Lemma for any given vector $v$ in a given highest weight $\hat{\mathfrak{g}}$-module $V$.
In Lemma \ref{lemma2} we take a Laurent polynomial $P(\zeta_1,\ldots,\zeta_s;\xi_1,\ldots,\xi_p)$ such that
$$R(\zeta,\xi)=\left(1+\zeta/\xi\right)^N$$
for some large positive integer $N$. Then for $r=0,1,\ldots,2s-1$ we have series (\ref{lemma2:1})
\begin{align}
&\lim_{\zeta,\xi\to\chi}D_{\zeta}^{r}\left(R(\zeta/\xi)X^{(s)}(\zeta)X^{(p)}(\xi)\right)\nonumber\\
=&\sum_{t=0}^{r}\binom{r}{t}\lim_{\zeta,\xi\to\chi}\left(D_{\zeta}^{t}\left(R(\zeta/\xi)\right)\left(D_{\zeta}^{r-t}X^{(s)}(\zeta)\right)X^{(p)}(\xi)\right).\label{lemma5:coeff}
\end{align}
Write
$$(j,j+1,\ldots,j+2s-1)=(j_1,j_2,\ldots,j_{2s}),$$
as in the proof of Lemma \ref{lemma3}.
By using Lemma \ref{lemma4} we can write a coefficient of (\ref{lemma5:coeff}) with summands (\ref{lemma5:summands}) as
\begin{align*}
\ldots& +\sum_{t=0}^{r}\binom{r}{t}\left(D_{\zeta}^{t}R\right)(1)j_{1}^{r-t} X^{(s)}(j_1) X^{(p)}(r-j_1)\\
&+\sum_{t=0}^{r}\binom{r}{t}\left(D_{\zeta}^{t}R\right)(1)j_{2}^{r-t} X^{(s)}(j_2) X^{(p)}(r-j_2)+\ldots\\
\ldots&+\sum_{t=0}^{r}\binom{r}{t}\left(D_{\zeta}^{t}R\right)(1)j_{2s}^{r-t} X^{(s)}(j_{2s}) X^{(p)}(r-j_{2s})+\ldots .
\end{align*}
Since
$$\sum_{t=0}^{r}\binom{r}{t}\left(D_{\zeta}^{t}R\right)(1)j^{r-t}=D_{z}^{r}\left(z^{j}R(z)\right)\Big|_{z=1}=a_{r,j}$$
for $j=j_1,\ldots,j_{2s}$ we can write a coefficient of (\ref{lemma5:coeff}) (i.e. a coefficient of (\ref{lemma2:1})) with summands
(\ref{lemma5:summands}) as
$$\ldots a_{r,1}X^{(s)}(j_1) X^{(p)}(r-j_1)+a_{r,2}X^{(s)}(j_2) X^{(p)}(r-j_2) +\ldots +a_{r,2s}X^{(s)}(j_{2s}) X^{(p)}(r-j_{2s})+\ldots .$$
By Lemma \ref{lemma3} the matrix 
$$\left(a_{r,m}\right)_{r=0,1,\ldots,2s-1;\hspace{2pt}m=j_1,j_2,\ldots,j_{2s}}$$
is regular, so we can express vectors (\ref{lemma5:summands}) in the way it is stated.
\end{prf}

The following corollary is an easy consequence of the above lemma.

\begin{cor}\label{cor6}
Let $n,s\in\mathbb{Z}$, $s>0$. For any given vector $v$ in any given highest weight $\hat{\mathfrak{g}}$-module $V$  every vector
$$X^{(s)}(j)X^{(s)}(n-j)v,\quad\textrm{where }n-j-2s<j\leq n-j,$$
in a formal Laurent series $X^{(s)}(\zeta)X^{(s)}(\xi)v$ can be written as a linear combination of vectors
$$X^{(s)}(i)X^{(s)}(n-i)v,\quad\textrm{where }i\leq n-i-2s,$$
and vectors obtained by the action of monomials, that have a coefficient of $X^{(s+1)}(\chi)$ as a factor, on $v$.
\end{cor}



\section{Quasi-particle bases of Verma modules}\label{5}
Let $V$ be a highest weight module with the highest weight $\Lambda$ and the highest weight vector $v_{\Lambda}$.
Recall the $\mathfrak{s}$-filtration (\ref{s:filtration}) and set
$$V_{(j)}=U_{(j)}v_{\Lambda}\subset V.$$
Then we have the {\em $\mathfrak{s}$-filtration} of the highest weight module $V$,
\begin{align}
&(0)=V_{(-1)}\subset V_{(0)}\subset V_{(1)}\subset \ldots \subset V,\label{s:filtration:module}\\
&V=\bigcup_{j\geq 0}V_{(j)}.\nonumber
\end{align}
Sometimes we shall write $$v\equiv w \textrm{ mod } V_{(j)}\quad\textrm{if}\quad v-w\in V_{(j)}.$$
Suppose $V$ is a Verma module $M=M(\Lambda)$ with the highest weight $\Lambda$ and the highest weight vector $v_{\Lambda}$. Then
in $M_{(p)}$, $p\geq 0$, we have a PBW basis:
\begin{align}
&\left\{B(i_1)\cdots B(i_r)X(j_1)\cdots X(j_s)v_{\Lambda} :\right.\nonumber\\
&\qquad \left. : r\geq 0,\, i_1\leq\ldots\leq i_r\leq -1,\, 0\leq s\leq p,\, j_1\leq\ldots\leq j_s\leq -1\right\}.\label{pbw:basis}
\end{align}
We shall say that
\begin{equation}\label{q-p:monomial}
X^{(p_1)}(j_1)\cdots X^{(p_s)}(j_s)
\end{equation}
is a {\em quasi-particle monomial}. Clearly
$$X^{(p_1)}(j_1)\cdots X^{(p_s)}(j_s)v_{\Lambda}\in M_{(p_1 +\ldots+p_s)}$$
and, as a consequence of (\ref{filtration:eq}), for any permutation $\sigma$
\begin{equation}\label{filtration:eq:charge}
X^{(p_{\sigma(1)})}(j_{\sigma(1)})\cdots X^{(p_{\sigma(s)})}(j_{\sigma(s)})v_{\Lambda}- X^{(p_1)}(j_1)\cdots X^{(p_s)}(j_s)v_{\Lambda}\in M_{(p_1 +\ldots+p_s-1)}.
\end{equation}
For this reason we consider quasi-particle monomials (\ref{q-p:monomial}) such that
\begin{align}
&p_1\leq \ldots\leq p_s,\label{cond:1}\\
&p_i =p_{i+1}\quad\textrm{implies}\quad j_i\leq j_{i+1},\label{cond:2}
\end{align}
and, by following Georgiev's terminology in \cite{Ge}, we say that such a quasi-particle monomial has
{\em charge-type} $$(p_1,\ldots,p_s),$$
{\em (total) charge} $$p_1+\ldots+p_s,$$
{\em degree-type} $$(j_1,\ldots,j_s)$$
and {\em (total) degree} $$j_1+\ldots+j_s.$$

\begin{lem}\label{verma:initial}
Let $p$ be a positive integer. Then
\begin{equation}\label{lemma:4:1}
X^{(p)}(j)v_{\Lambda}\in M_{(p-1)}\quad\textrm{for}\quad j>-p.
\end{equation}
\end{lem}

\begin{prf}
Let $p=1$. Then $X^{(1)}(0)=X(0)=h\otimes t^0$ and $X(0)v_{\Lambda}=\Lambda(h)v_{\Lambda}\in M_{(0)}$.

Now assume that (\ref{lemma:4:1}) holds for $p$. Recall recursive formula (\ref{Xp:recursive}).
Then for given $j>-p-1$, by using (\ref{filtration:eq:charge}), we get (for some coefficients $\alpha_i$) a finite sum 
\begin{align*}
X^{(p+1)}(j)v_{\Lambda}&\equiv\sum_{j-i\leq -p}\alpha_{i}X(i)X^{(p)}(j-i)v_{\Lambda}\quad\textrm{mod }M_{(p)}\\
&\equiv \sum_{\substack{j-i\leq -p\\i\leq -1}}\alpha_{i}X^{(p)}(j-i)X(i)v_{\Lambda}\quad\textrm{mod }M_{(p)}\\
&\equiv 0\quad\textrm{mod }M_{(p)}.
\end{align*}
\end{prf}

As in \cite{Ge} we introduce a linear order ``$\prec$'' on the set of monomials (\ref{q-p:monomial}) of quasi-particles (satisfying (\ref{cond:1}) and (\ref{cond:2})). 
Let $b$ and $b'$ be two quasi-particle monomials. The order ``$\prec$'' is defined as follows:
We shall write $b\prec b'$ 
\renewcommand{\labelitemi}{$\circ$}
\begin{itemize}
  \item if total charge of $b$ is greater than total charge of $b'$;
  \item if total charges of $b$ and $b'$ are the same and charge-type of $b$ is less than charge-type of $b'$ by reverse lexicographic order;
  \item if charge-types of $b$ and $b'$ are the same and degree-type of $b$ is less then degree-type of $b'$ by reverse lexicographic order.
\end{itemize}

\begin{ex}
We have
$$X^{(2)}(j_1)X^{(3)}(j_2)\prec X^{(1)}(i_1)X^{(2)}(i_2)$$
since $2+3>1+2$;
$$X^{(2)}(j_1)X^{(3)}(j_2)\prec X^{(1)}(i_1)X^{(4)}(i_2)$$
since $2+3=1+4$ and $3<4$; and
$$X^{(2)}(-2)X^{(2)}(-2)X^{(3)}(-5)\prec X^{(2)}(-3)X^{(2)}(-1)X^{(3)}(-5)$$
since $(2,2,3)=(2,2,3)$, $-5=-5$ and $-2<-1$.
\end{ex}

The main result in this section is

\begin{thm}\label{verma:main}
The set $\mathcal{B}_{M(\Lambda)}$  of vectors
\begin{equation}\tag{V}\label{verma:basis}
B(i_1)\cdots B(i_r)X^{(p_1)}(j_1)\cdots X^{(p_s)}(j_s)v_{\Lambda}
\end{equation}
such that
\begin{align}
 &r\geq 0 \quad\textrm{and odd}\quad i_1\leq \ldots\leq i_r\leq -1,\tag{V1}\label{v1}\\
&s\geq 0\quad\textrm{and}\quad 1\leq p_1\leq \ldots\leq p_s,\tag{V2}\label{v2}\\
& j_s\leq -p_s,\tag{V3}\label{v3}\\
&p_l<p_{l+1}\quad\textrm{implies}\quad j_l\leq -p_l -2p_{l}(s-l),\tag{V4}\label{v4}\\
&p_l=p_{l+1}\quad\textrm{implies}\quad j_l\leq -2p_l + j_{l+1}\tag{V5}\label{v5}
\end{align}
is a basis of the Verma module $M(\Lambda)$.
\end{thm}

\begin{prf}
First we prove that $\mathcal{B}_{M(\Lambda)}$ is a spanning set for $M(\Lambda)$. We use the induction on our order ``$\prec$''. 
We start with a PBW spanning set of vectors
\begin{align}
&B(i_1)\cdots B(i_r)X^{(1)}(j_1)\cdots X^{(1)}(j_s)v_{\Lambda},\label{spanning:set}\\
 &r\geq 0 \quad\textrm{and}\quad i_1\leq \ldots\leq i_r\leq -1,\nonumber\\
&s\geq 0\quad\textrm{and}\quad j_1\leq \ldots\leq j_s\leq -1,\nonumber
\end{align}
If for a vector (\ref{spanning:set}) the difference condition (\ref{v5}), $j_l\leq -2+j_{l+1}$,
is not satisfied for all $l=1,2,\ldots,s-1$, we use Corollary \ref{cor6} to
replace $X^{(1)}(j)X^{(1)}(j)$ or $X^{(1)}(j-1)X^{(1)}(j)$ with a linear combination
of monomials $X^{(1)}(i_1)X^{(1)}(i_2)$, $i_1\leq -2+i_2$, and a quasi-particle $X^{(2)}(i)$.
As a result we shall replace a vector for which the difference condition is not satisfied by quasi-particle
monomials  which are higher in order ``$\prec$''. 
In fact, by using Lemma \ref{lemma5}, Corollary \ref{cor6} and Lemma \ref{verma:initial}, and by repeating
Georgiev's argument in \cite{Ge} (almost) verbatim, we see that the set of vectors $\mathcal{B}_{M(\Lambda)}$ is a spanning 
set for $M(\Lambda)$.

In order to prove linear independence of the set of vectors $\mathcal{B}_{M(\Lambda)}$ it is sufficient to see that
for every integer $n$ the number of vectors $\mathcal{B}_{M(\Lambda)}$ of degree $n$ is equal to dimension
$\dim M_{n}$ given by  character formula (\ref{verma:char}). But this statement follows from Georgiev's construction
of quasi-particle bases of principal subspaces of standard $\widehat{\mathfrak{sl}}_{2}$-modules
$L(k\Lambda_{0})$ in the homogeneous picture---it can be stated as a character formula  (cf. \cite{Bu} and Section \ref{7}):
$$\ch\textstyle_{q} M(\Lambda) = \displaystyle F\cdot \sum_{n_1,n_2,\ldots} \frac{q^{N_{1}^{2}+N_{2}^{2}+\ldots} }{ (q)_{n_1} (q)_{n_2}\cdots },$$
where $N_j = n_j +n_{j+1}+\ldots$ and the sum goes over all sequences $n_1,n_2,\ldots$ of nonnegative integers with finitely many nonzero entries.
\end{prf}

We introduce a linear order ``$\prec$'' on the set of basis vectors $\mathcal{B}_{M(\Lambda)}$ given by Theorem 
\ref{verma:main}: We first compare two basis vectors by comparing their quasi-particle monomial factors, and if 
equal, we then compare monomials in Heisenberg Lie algebra elements $B(j)$ by reverse lexicographical order.

If vector $v\in M(\Lambda)$, $v\neq 0$, is written as a linear combination
of basis elements $\mathcal{B}_{M(\Lambda)}$ of $M(\Lambda)$,
\begin{equation}\label{pokemonlopta}
v=\lambda_{b}b+\sum_{b\prec b'} \lambda_{b'}b',\quad \lambda_{b}\neq 0,
\end{equation}
then we say that $b$ is the {\em leading term} of $v$ and write $b=\lt (v)$.


\section{Quasi-particle bases of standard modules}

Let $E^{+}(\zeta)$ and $E^{-}(\zeta)$ be the formal Laurent series with coefficients in $U=U(\hat{\mathfrak{g}})$ defined as the formal exponential series
\begin{equation*}
E^{\pm}(\zeta)=\sum_{i>0} E^{\pm}(\pm i)\zeta^{\pm i}=\exp\left(2\sum_{n\in\pm(2\mathbb{N}+1)}B(n)\zeta^{n}/n\right).
\end{equation*}
Let $k_0$, $k_1$ be nonnegative integers and set $k=k_0+k_1$. Then on the standard $\hat{\mathfrak{g}}$-module $L(k_0\Lambda_0 + k_1\Lambda_1)$
we have (see \cite[Theorem 6.6.]{MP}):
\begin{itemize}
  \item [(a)] For $p,q\geq 0$, $p+q=k$,
  \begin{equation}\label{a}
  a_p X^{(p)}(\zeta)-(-1)^{k_0}a_q E^{-}(-\zeta)X^{(q)}(-\zeta)E^{+}(-\zeta)=0,
  \end{equation}
  where $a_r=2^{-r(r-2)}/r!$.
  \item [(b)] For $p\geq k+1$,
  \begin{equation}\label{b}
  X^{(p)}(\zeta)=0.
  \end{equation}
\end{itemize}

Let $v_{\Lambda}$ be a highest weight vector in $L=L(\Lambda)$. Then we have the $\mathfrak{s}$-filtration (recall (\ref{s:filtration:module})),
$$L_{(j)}=U_{(j)}v_{\Lambda}\subset L \qquad \textrm{for }j\in\mathbb{Z}.$$

\begin{lem}\label{lemma:5:1}
Set $t=\min\left\{k_0,k_1\right\}$. Then
\begin{itemize}
  \item [(a)] for $r=1,2,\ldots,t$ we have
  \begin{equation*}
  X^{(r)}(l)v_{\Lambda}\in L_{(r-1)}\qquad\textrm{for }l>-r,
  \end{equation*}
  \item [(b)] for $r=t+1,\ldots,[k/2]$ we have
  \begin{equation*}
  X^{(r)}(l)v_{\Lambda}\in L_{(r-1)}\qquad\textrm{for }l>-2r+t.
  \end{equation*}
\end{itemize}
\end{lem}

\begin{prf}
(a) Since $L(\Lambda)$ is a quotient of $M(\Lambda)$, the statement follows from Lemma \ref{verma:initial}.

\noindent (b) From the proof of  Proposition 8.6. in \cite{MP} we have
\begin{equation}\label{ł}
X(-1)^{t+1}v_{\Lambda}\in L_{(t)}.
\end{equation}

Let $r=t+s$, $1\leq s\leq [k/2]-t$, $l>-t-2s$. Then we have a finite sum (with some coefficients $a_{j_1,\ldots,j_{t+s}}$)
$$X^{(t+s)}(l)v_{\Lambda}=\sum_{j_1+\ldots+j_{t+s}=l} a_{j_1,\ldots,j_{t+s}}X(j_1)\cdots X(j_{t+s})v_{\Lambda}.$$
By using (\ref{filtration:eq}) we see that
$$X^{(t+s)}(l)v_{\Lambda}\equiv\sum_{\substack{j_1+\ldots+j_{t+s}=l\\j_1,\ldots,j_{t+s}\leq -1}} a_{j_1,\ldots,j_{t+s}}X(j_1)\cdots X(j_{t+s})v_{\Lambda}\quad\textrm{mod }L_{(r-1)}.$$
Now we have $-t-2s<l\leq -t-s$, so at least $t+1$ of indices $j_1,\ldots,j_{t+s}$ must be equal to $-1$. This means that, modulo $U_{(r-1)}$,
each product $X(j_1)\cdots X(j_{t+s})$ contains a factor $X(-1)^{t+1}$, and (b) follows from (\ref{ł}).
\end{prf}

Set  $$t=\min\left\{k_0,k_1\right\}$$ and for $p\leq [k/2]$ set
\begin{align*}n_{\Lambda}(p)=\left\{\begin{array}{l@{\,\ }l}
 p &\textrm{ for }p\leq t,\\
 2p-t&\textrm{ for }p\geq t+1.
\end{array}\right.
\end{align*}

The main result in this paper is
\begin{thm}\label{standard:main}
The set $\mathcal{B}_{L(\Lambda)}$ of vectors
\begin{equation}\tag{S}\label{standard:basis}
B(i_1)\cdots B(i_r)X^{(p_1)}(j_1)\cdots X^{(p_s)}(j_s)v_{\Lambda}
\end{equation}
such that
\begin{align}
 &r\geq 0 \quad\textrm{and odd}\quad i_1\leq \ldots\leq i_r\leq -1,\tag{S1}\label{s1}\\
&s\geq 0\quad\textrm{and}\quad 1\leq p_1\leq \ldots\leq p_s\leq [k/2],\tag{S2}\label{s2}\\
& j_s\leq -n_{\Lambda}(p_s),\tag{S3}\label{s3}\\
&p_l<p_{l+1}\quad\textrm{implies}\quad j_l\leq -n_{\Lambda}(p_l) -2p_{l}(s-l),\tag{S4}\label{s4}\\
&p_l=p_{l+1}\quad\textrm{implies}\quad j_l\leq -2p_l + j_{l+1}\tag{S5}\label{s5}\\
&\textrm{if }k\textrm{ is even and }p_l=k/2,\textrm{ then }j_l\equiv k_0(\textrm{mod }2)\tag{S6}\label{s6}
\end{align}
is a basis of the standard module $L(k_0 \Lambda_{0}+k_1\Lambda_1)$.
\end{thm}

\begin{prf}
Since $L(\Lambda)$ is a quotient of $M(\Lambda)$, the set of vectors satisfying (\ref{v1})--(\ref{v5}) spans $L(\Lambda)$.

While for Verma module all charges $p$ for quasi-particle $X^{(p)}(j)$ are allowed in Theorem \ref{verma:main}, here we have a restriction in (\ref{s2}):
$$p_s\leq [k/2].$$
This is clear for $p\geq k+1$ since by (\ref{b}) $X^{(p)}(j)=0$ on $L$. 
For $k/2<p\leq k$ we use (\ref{a}) to express $X^{(p)}(j)$ in terms of operators $X^{(q)}(i)$ which are lower in $\mathfrak{s}$-filtration
(since $q=k-p<k/2$). 

Finally, for $k$ even and $p=k/2$ relation (\ref{b}) can be written as (see \cite[relation (8.3)]{MP})
\begin{align*}
2X^{(p)}(n)&=\left(1-(-1)^{k_0 +n}\right)X^{(p)}(n)\\
&=-\sum_{i<0}E^{-}(i)X^{(p)}(n-i)-\sum_{i>0}X^{(p)}(n-i)E^{+}(i)
\end{align*}
for $n\not\equiv k_0(\textrm{mod }2)$. Therefore, we can replace a vector which contains $X^{(p)}(n)$ with other vectors higher in our order.

For $k_0=[k/2]$ we have
$$t=\min\left\{[k/2],k-[k/2]\right\}=[k/2]\quad\textrm{and}\quad n_{\Lambda}(p)=p$$
so conditions (\ref{v4}) and (\ref{s4}) in Theorems \ref{verma:main} and \ref{standard:main} coincide. 
So in this case we see that vectors (\ref{standard:basis}) span $L(\Lambda)$. In other cases we have to repeat Georgiev's argument:
We consider the case
$$p_l<p_{l+1}$$
and conclude, by using (\ref{filtration:eq:charge}) and Lemma \ref{lemma:5:1}, that a monomial containing
$$X^{(p_l)}(j_l),\qquad j_{l}>-n_{\Lambda}(p_l)$$
can be replaced with monomials higher in our order. Then we consider all the factors ``in front of''
$X^{(p_l)}(j_l)$:
$$\cdots X^{(p_l)}(j_l)X^{(p_{l+1})}(j_{p_{l+1}})\cdots X^{(p_s)}(j_s)v_{\Lambda}.$$
By using (\ref{filtration:eq:charge}) and Lemma \ref{lemma5} we conclude that a monomial containing
$$X^{(p_l)}(j_l),\qquad j_l>-n_{\Lambda}(p_l)-2p_l$$
can be replaced with monomials higher in our order. After repeating the same argument
for all the factors in front of
$X^{(p_l)}(j_l)$, we obtain the condition (\ref{s4}).

Hence the set of vectors (\ref{standard:basis}) satisfying conditions (\ref{s1})--(\ref{s6}) is a spanning set of $L(\Lambda)$.
In the next section we show that the number of vectors of degree $n$ in the spanning set is given by $n$th coefficient of the sum side of
(\ref{grr}), while the product side is the principally specialized character of $L(\Lambda)$. Hence linear independence follows.
\end{prf}


\section{The Rogers-Ramanujan-type identities}\label{7}

Let  $L=L(\Lambda)$ be a standard $\widehat{\mathfrak{sl}}_{2}$-module with the highest weight
$$\Lambda=k_0\Lambda_0+ k_1\Lambda_1$$
and set $$k=k_0+k_1.$$
  We have a principally specialized character formula, that follows from the Weyl-Kac formula (see \cite{LM}),
\begin{align*}\ch\textstyle_{q} L(\Lambda)=\left\{\begin{array}{l@{\,\ }l}
 \displaystyle F\cdot\prod_{\substack{n\geq 1\\n\not\equiv 0,\pm (k_0 +1)}}(1-q^n)^{-1} &\textrm{ if }k_0\neq k_1,\\
 \displaystyle F\cdot\prod_{\substack{n\geq 1\\n\not\equiv 0, (k_0 +1)}}(1-q^n)^{-1} \prod_{\substack{n\geq 1\\n\equiv k_0 +1}}(1-q^n)\hspace{10pt}&\textrm{ if }k_0=k_1,
\end{array}\right.
\end{align*}
where all congruences are modulo $k+2$ and
$$F=\prod_{n\geq 1}(1-q^{2n-1})^{-1}.$$ 

For nonnegative integers $n_1,\ldots,n_r$ let
$$N_{j}=n_j + n_{j+1} +\ldots+n_r.$$
Set $$i=\min\left\{k_0 +1,k_1+1\right\}.$$ 
If $k$ is odd, we have Andrews' generalization of Rogers-Ramanujan identities   (cf. \cite[Theorem 7.8.]{A2}),
\begin{equation}\label{andrews}
\prod_{\substack{n\geq 1\\n\not\equiv 0,\pm (k_0 +1)}}(1-q^n)^{-1} = 
\sum_{n_1,n_2,\ldots,n_{[k/2]}\geq 0} \frac{q^{N_{1}^{2}+N_{2}^{2}+...+ N^{2}_{[k/2]}+N_{i}+N_{i+1}+\ldots+N_{[k/2]}   }}{(q)_{n_1} (q)_{n_2}\cdots (q)_{n_{[k/2]}}}.
\end{equation}
If $k$ is even and $k_0\neq k_1$ we have Bressoud's generalization of Rogers-Ramanujan identities (cf. \cite{Br2}),
\begin{equation}\label{bressoud}
\prod_{\substack{n\geq 1\\n\not\equiv 0,\pm (k_0 +1)}}(1-q^n)^{-1} = 
\sum_{n_1,n_2,\ldots,n_{k/2}\geq 0} \frac{q^{N_{1}^{2}+N_{2}^{2}+...+ N^{2}_{k/2}+N_{i}+N_{i+1}+\ldots+N_{k/2}   }}{(q)_{n_1} (q)_{n_2}\cdots (q)_{n_{k/2-1}}   (q^2)_{n_{k/2}}}.
\end{equation}
If $k_0=k_1$ we have another generalization obtained by Bressoud (cf. \cite{Br1}),
\begin{equation}\label{bressoud2}
\prod_{\substack{n\geq 1\\n\not\equiv 0, (k_0 +1)}}(1-q^n)^{-1} \prod_{\substack{n\geq 1\\n\equiv k_0 +1}}(1-q^n) = 
\sum_{n_1,n_2,\ldots,n_{k/2}\geq 0} \frac{q^{N_{1}^{2}+N_{2}^{2}+...+ N^{2}_{k/2}   }}{(q)_{n_1} (q)_{n_2}\cdots (q)_{n_{k/2-1}}   (q^2)_{n_{k/2}}}.
\end{equation}

The product of right sides of identities (\ref{andrews}), (\ref{bressoud}) and (\ref{bressoud2}) with $F=\prod_{n\geq 1}(1-q^{2n-1})^{-1}$ 
corresponds to the bases of the 
standard modules $L(\Lambda)$, given by Theorem \ref{standard:main}, in a way that each $n_j$ 
on the right side
corresponds to the quasi-particle of charge $j$.
Therefore, these identities can be used in order to immediately obtain linear independence of the set (\ref{standard:basis}).
Likewise, they can be considered as corollaries of  Theorem \ref{standard:main}, if we give another proof of linear independence of the set (\ref{standard:basis}). 
In the rest of this section we will briefly explain the above-mentioned correspondence between quasi-particle charges and degrees and the sum side of identities 
(\ref{andrews}), (\ref{bressoud}) and (\ref{bressoud2}).

The number of monomials 
$$B(i_1)\cdots B(i_r),\quad i_1+\ldots+i_r=m,$$ 
in (\ref{standard:basis}), satisfying condition (\ref{s1}),
$$r\geq 0 \quad\textrm{and odd}\quad i_1\leq \ldots\leq i_r\leq -1,$$
equals to the $n$th coefficient in 
$$F=\prod_{n\geq 1}(1-q^{2n-1})^{-1}=\prod_{n\geq 1}(1+q^{2n-1}+q^{2(2n-1)}+q^{3(2n-1)}+\ldots).$$ 
In the argument above we implicitly identify monomials and partitions:
$$B(-j_1)\cdots B(-j_r)\qquad \longleftrightarrow\qquad j_1+\ldots+j_r,\textrm{ }j_1\geq \ldots\geq j_r\geq 1.$$

Let us start with $\Lambda=k\Lambda_0$, i.e. $k_0=k$, $k_1=0$, $n_{\Lambda}(p)=2p$. Let $k$ be odd.
For quasi-particle monomials of charge-type
$$(n_1,n_2,\ldots,n_{[k/2]})$$
term
\begin{equation}\label{RR_term}
\frac{q^{2p\cdot \tfrac{1}{2}n_p (n_p -1) +n_p 2p + n_p 2p(n_{p+1}+\ldots+n_{[k/2]})   }}{(q)_{n_p}}
\end{equation}
``counts'' products of $n_p$ quasi-particles of charge $p$ satisfying difference conditions (\ref{s5})
and initial conditions (\ref{s4}). Note that for $p=[k/2]$ we have the term
$$\frac{q^{pn_{p}^{2}+pn_{p}}}{(q)_{n_p}}.$$
The product of all terms (\ref{RR_term}) for $p=1,2,\ldots,[k/2]$ gives a
summand in (\ref{andrews}) for $i=1$.
In the case when $k$ is even we have $[k/2]=k/2$ and for $p=k/2$ in (\ref{standard:basis})
we have only quasi-particles
$$X^{(p)}(j),\quad j\equiv k\textrm{ mod }2$$
since $k_0 =k$.
Note that the ``smallest'' partition satisfying difference and initial condition is
\begin{equation}\label{RR_partition}
2p,\, 2p+2p,\, \ldots,\, 2p+2p(n_{p}-1)
\end{equation}
so all the parts are congruent  $2\cdot\tfrac{k}{2}$ \textrm{mod} $2$,
and we obtain all the other by ``adding'' partitions with even parts only, represented by $1/(q^2)_{n_p}$.
So for $k$ even the last term is
$$\frac{q^{pn_{p}^{2}+pn_{p}}}{(q^2)_{n_p}},\quad p=k/2,$$
and the product of all terms gives the summand in (\ref{bressoud2}) for $i=1$.

Recall that for $\Lambda =k\Lambda_{0}$ we have $i=1$ and initial conditions 
$$X^{(p)}(j),\quad j\leq -2p.$$

In the case $\Lambda=(k-1)\Lambda_{0}+\Lambda_{1}$ we have $i=2$ and initial conditions
$$X^{(1)}(j),\quad j\leq -1=-2+1,\qquad X^{(p)}(j),\quad j\leq -2p+1\textrm{ for }p\geq 2.$$
This means that for each $p=1,\ldots,[k/2]$ the parts of the ``smallest'' partition satisfying 
difference conditions and initial conditions are smaller for $1$, i.e., they are
\begin{equation}\label{RR_partition2}
2p-1,\, 2p-1+2p,\, \ldots,\, 2p-1+2p(n_{p}-1)
\end{equation}
(in fact, each part is shifted by $2p(n_{p+1}+\ldots+n_{[k/2]})$).
This gives for $i=2$ a term ``smaller for'' $N_1=n_1+n_2+\ldots+n_{[k/2]}$,
$$q^{N_{1}^{2}+...+ N^{2}_{[k/2]}+N_{2}+\ldots+N_{[k/2]} },$$
when compared with the previous term
$$q^{N_{1}^{2}+...+ N^{2}_{[k/2]}+N_1 +N_{2}+\ldots+N_{[k/2]} }$$
for $i=1$.
Also note that, for $k$ even and $p=k/2$, all parts in (\ref{RR_partition2}) are odd, just like $k_0 =k-1$.

In the case $\Lambda=(k-2)\Lambda_{0}+2\Lambda_1$ we have $i=3$ and initial conditions
$$X^{(1)}(j),\quad j\leq -1,\qquad X^{(2)}(j),\quad j\leq -2=-4+2,\qquad X^{(p)}(j),\quad j\leq -2p+2\textrm{ for }p\geq 3.$$
This means that for each $p=2,\ldots,[k/2]$ the parts of the ``smallest'' partition satisfying difference conditions 
and initial conditions are smaller for $1$ when compared with the previous case $\Lambda=(k-1)\Lambda_{0}+\Lambda_1$
and $i=2$.
This gives for $i=3$ a term ``smaller for''
$N_2=n_2+\ldots + n_{[k/2]}$,
$$q^{N_{1}^{2}+...+ N^{2}_{[k/2]}+N_{3}+\ldots+N_{[k/2]} },$$
when compared with the previous term
$$q^{N_{1}^{2}+...+ N^{2}_{[k/2]}+N_{2}+N_{3}+\ldots+N_{[k/2]} }$$
for $i=2$.

In a similar way we see that for all possible $\Lambda=k_0\Lambda_{0}+k_1\Lambda_{1}$
summands in (\ref{andrews})--(\ref{bressoud2}) count quasi-particle
monomials of charge-type $(n_1,n_2,\ldots,n_{[k/2]})$
satisfying (\ref{s3})--(\ref{s6}).

\section{Bases of maximal submodules of some Verma modules}\label{sec:last}

As in \cite{MP}, we would like to prove linear independence of the basis of $L(\Lambda)$
by vertex-operator theoretic methods, i.e., by constructing a ``complementary'' basis of the maximal submodule $W(\Lambda)$
of Verma module $M(\Lambda)$. Unfortunately, we do not understand well the role of initial conditions in the structure of $W(\Lambda)$,
and we are able to construct desired bases only for odd level modules
$W((l+1)\Lambda_{0}+l\Lambda_{1})$ and $W(l\Lambda_{0}+(l+1)\Lambda_{1})$, $l\geq 1$.

Set $a_p =0$ for $p<0$.
As in \cite{MP}, we define $R_{p}^{\Lambda}(n)$ as coefficients in
$$R_{p}^{\Lambda}(\zeta)=\sum_{n\in\mathbb{Z}}R_{p}^{\Lambda}(n)\zeta^{n},$$
where
\begin{align*}
R_{p}^{\Lambda}(\zeta)=a_p X^{(p)}(\zeta) - (-1)^{k_0}a_q E^{-}(-\zeta)X^{(q)}(-\zeta)E^{+}(-\zeta),\quad p+q=k.
\end{align*}
For $p,n\in\mathbb{Z}$ we also write
$$R_{p}(n)=R_{p}^{\Lambda}(n),\quad R_{p}(\zeta)=R_{p}^{\Lambda}(\zeta).$$
By (\ref{a}) and (\ref{b}) we have $R_{p}^{\Lambda}(n)=0$ on $L(\Lambda)$, so for a highest weight vector $v_{\Lambda}\in M(\Lambda)$ we have
$$R_{p}^{\Lambda}(n)v_{\Lambda}\in W(\Lambda)\quad \textrm{for all }p\in\mathbb{Z},\, n\geq 0.$$

\begin{lem}
(a) For $p>k/2$ and $j\leq -p$ $$\lt (R_{p}^{\Lambda}(j)v_{\Lambda})=X^{(p)}(j)v_{\Lambda}.$$
(b) If $k$ is even, $p=k/2$ and $j\leq -p$, $j\not\equiv k_0\,(\text{mod }2)$
$$\lt(R_{p}^{\Lambda}(j)v_{\Lambda})=X^{(p)}(j)v_{\Lambda}. $$
\end{lem}

\begin{prf}
(a) If $p\geq k+1$, we have by definition
$$R_{p}^{\Lambda}(j)v_{\Lambda}=a_p X^{(p)}(j)v_{\Lambda}$$
and the statement is clear.
If $p\leq k$ and $p+q=k$, then we have $p>k/2>q$ and the coefficient of $\zeta ^j$ in
$E^{-}(-\zeta)X^{(q)}(-\zeta)E^{+}(-\zeta)v_{\Lambda}$ is a linear combination of basis elements in 
$M_{(q)}$ and hence higher in order than $X^{(p)}(j)$.

(b)  If $k$ is even and $p=k/2$, then
$$R_{p}^{\Lambda}(j)v_{\Lambda}=a_{p}\left(X^{(p)}(j)-(-1)^{k_0 +j}X^{(p)}(j)-(-1)^{k_0 +j}\sum_{i>0}E^{-}(-i)X^{(p)}(j+i)\right)v_{\Lambda}$$
and $X^{(p)}(j)v_{\Lambda}$ is the leading term when $k_0 +j$ is odd. 
\end{prf}

From Lemma \ref{lemma:5:1} we see that for $p=1,\ldots,[k/2]$ there exist elements
\begin{equation}\label{shorter}
I_{p}^{\Lambda}(n)\equiv X^{(p)}(n)\quad\textrm{modulo ``shorter'' terms in $\mathfrak{s}$-filtration}
\end{equation}
such that on $L(\Lambda)$
$$I_{p}^{\Lambda}(n)v_{\Lambda}=0\qquad\textrm{for }n>-n_{\Lambda}(p).$$
Hence for a highest weight vector $v_{\Lambda}\in M(\Lambda)$ and $t+1\leq p\leq [k/2]$ we have
$$I_{p}^{\Lambda}(n)v_{\Lambda}\in W(\Lambda)\quad\text{for }-2p+t<n\leq -p.$$
It is clear from the construction, i.e., from (\ref{shorter}), that 
$$\lt(I_{p}^{\Lambda}(n)v_{\Lambda})=X^{(p)}(n)v_{\Lambda}\quad \text{for }-2p+t<n\leq -p.$$

Consider the set of vectors 
\begin{equation}\label{setofvectors}
v=B(i_1)\cdots B(i_r)X^{(p_1)}(j_1)\cdots X^{(p_s)}(j_s)D_{p}^{\Lambda}(n)v_{\Lambda},
\end{equation}
where $D_{p}^{\Lambda}(n)$ is either
$$R_p^{\Lambda}(n),\quad p\geq [k/2],\, n\leq -p,$$
with $j+k_0$ odd if $p=k/2$, or
$$I_{p}^{\Lambda}(n),\quad t+1\leq p\leq [k/2],$$
such that
$$\lt v\equiv B(i_1)\cdots B(i_r)X^{(p_1)}(j_1)\cdots X^{(p_s)}(j_s)X^{(p)}(n)v_{\Lambda}\quad\text{mod }M_{(p_1+\ldots+p_s +p-1)}.$$
Our assumption implies that the quasi-particle monomial vector
\begin{equation}\label{led_term}
\lt v= B(i_1)\cdots B(i_r)X^{(p_1)}(j_1)\cdots X^{(p)}(n)\cdots  X^{(p_s)}(j_s)v_{\Lambda},\quad p_1\leq\ldots\leq p\leq\ldots\leq p_s,
\end{equation}
satisfies difference and initial conditions (\ref{v1})--(\ref{v5}).

Choose any set of vectors $\mathcal{B}_{W(\Lambda)}$ of the form (\ref{setofvectors}) such that each leading term (\ref{led_term})
appears only once. Then (\ref{pokemonlopta}) implies that $\mathcal{B}_{W(\Lambda)}$ is a linearly independent set.
The set $\mathcal{B}_{W(\Lambda)}$ is ``complementary'' to the basis $\mathcal{B}_{L(\Lambda)}$ in a sense that
\begin{equation}\label{srce}
\mathcal{B}_{M(\Lambda)}=\mathcal{B}_{L(\Lambda)}' \cup \lt\mathcal{B}_{W(\Lambda)},\quad \mathcal{B}_{L(\Lambda)}' \cap \lt\mathcal{B}_{W(\Lambda)}=\emptyset,
\end{equation}
for a set of basis vectors $\mathcal{B}_{L(\Lambda)}'$ satisfying difference and initial conditions (\ref{s1})--(\ref{s6}).
Then Theorem \ref{standard:main} and (\ref{srce}) imply that
\begin{equation}\label{tref}
\mathcal{B}_{W(\Lambda)}\text{ is a basis of }W(\Lambda).
\end{equation}
Note that (\ref{tref}) is equivalent to the statement
\begin{equation}\label{doubletref}
\lt\mathcal{B}_{W(\Lambda)}=\left\{\lt v : v\in W(\Lambda),\, v\neq 0\right\}.
\end{equation}
As we have already said, we would like to construct a ``complementary'' basis $\mathcal{B}_{W(\Lambda)}$ of $W(\Lambda)$ by
vertex-operator theoretic methods---in the setting above it would be enough to prove that $\mathcal{B}_{W(\Lambda)}$
is a spanning set for $W(\Lambda)$.

As a consequence of Proposition 9.1 in \cite{MP} we have
\begin{pro}\label{8.1}
For $p,q\in\mathbb{Z}$, $p+q=k$,
$$R_p (\zeta)v_{\Lambda}=-(-1)^{k_0}E^{-}(-\zeta)R_{q}(-\zeta)v_{\Lambda}.$$
\end{pro}

By Proposition 9.2 in \cite{MP} for $p\in\mathbb{Z}$ and $p+q=k$ on a highest weight module we have
\begin{align}
[X(\zeta),R_{p}(\chi)]=&2(q+1)D\delta(-\zeta/\chi)R_{p-1}(\chi)\nonumber\\
&-4\delta(-\zeta/\chi)\left(B^{-}(\chi)R_{p-1}(\chi)+R_{p-1}(\chi)B^{+}(\chi)\right)\label{1:1}\\
&-2\delta(-\zeta/\chi)D_{\chi}R_{p-1}(\chi).\nonumber
\end{align}

By Proposition 9.3 in \cite{MP} we have
\begin{align}
\left(\frac{1+\zeta/\chi}{1-\zeta/\chi}\right)^{2} X(\zeta)R_{p}(\chi)-\left(\frac{1+\chi/\zeta}{1-\chi/\zeta}\right)^{2}R_{p}(\chi)X(\zeta)\nonumber\\
=2(p+1)D\delta(\zeta/\chi)R_{p+1}(\chi)-2\delta(\zeta/\chi)D_{\chi}R_{p+1}(\chi).\label{2:2}
\end{align}
Recall the Taylor expansion 
\begin{equation*}
\left(\frac{1+z}{1-z}\right)^2 = 1+\sum_{n\geq 1}4n z^n .
\end{equation*}

\begin{pro}\label{8.2}
Let $p,M\in\mathbb{Z}$. Then
\begin{align*}
-2MR_{p+1}(M)v_{\Lambda}=&-2MR_{p-1}(M)v_{\Lambda}-\sum_{\text{odd }i>0}B(-i)R_{p-1}(M+i)v_{\Lambda}\\& +\sum_{r\geq 1}4rX(-r)R_{p}(M+r)v_{\Lambda}.
\end{align*}
\end{pro}

\begin{prf}
By (\ref{2:2}) we have
\begin{align*}
[X(\zeta),R_{p}(\chi)&]+\left(\sum_{r\geq 1}4r\left(\zeta/\chi\right)^{r}\right)X(\zeta) R_{p}(\chi) - \left(\sum_{r\geq 1}4r\left(\chi/\zeta\right)^{r}\right)R_{p}(\chi)X(\zeta)\\
&=2(p+1)D\delta(\zeta/\chi)R_{p+1}(\chi)-2\delta(\zeta/\chi)D_{\chi}R_{p+1}(\chi).
\end{align*}
We use (\ref{1:1}) to express $[X(\zeta),R_{p}(\chi)]$ in terms of $R_{p-1}(\chi)$ and then equate the coefficients of $\zeta^0 \chi^M$.
\end{prf}

By Proposition 9.4 in \cite{MP} we have
\begin{equation}\label{3:3}
[B(j),R_{p}(n)]=2pR_{p}(j+n)
\end{equation}
for all $j,p,n\in\mathbb{Z}$.
Set
$$R^{\Lambda}=U(\mathfrak{s}_{-})\spn\left\{R_{p}^{\Lambda}(n)v_{\Lambda} : p,n\in\mathbb{Z}\right\}.$$
Note that the maximal submodule $W(\Lambda)$ is generated by two singular vectors (cf. \cite{K})
$$f_{0}^{k_0 +1}v_{\Lambda}\quad\text{and}\quad f_{1}^{k_1 +1}v_{\Lambda}$$
of degrees $-k_0 -1$ and $-k_1 - 1$, so
$$R_{p}^{\Lambda}(n)v_{\Lambda}=0\quad \text{for }n\geq -\min\left\{k_0 ,k_1\right\}.$$
The commutation relations (\ref{1:1}) and (\ref{3:3}) imply that 
$$\mathfrak{n}_{+}R^{\Lambda}\subset R^{\Lambda},$$
so
$$V''=U(\mathfrak{g})R^{\Lambda}=U(\mathfrak{n}_{-})R^{\Lambda}\subset W(\Lambda)$$
is $\widehat{sl}_2$-submodule of $W(\Lambda)$.
In fact, it is proved in \cite{MP} that $V''=W(\Lambda)$ by using a construction of another combinatorial basis 
of $V''$. As we have already announced, in this section we'll construct  bases of some $W(\Lambda)$.
So from now on we fix
\begin{equation}\label{lambda}
\Lambda=(l+1)\Lambda_0 + l\Lambda_1,\quad k=2l+1,\, l\geq 1.
\end{equation}
Then
$$R_{p}^{\Lambda}(n)v_{\Lambda}=0\quad \textrm{for }n\geq -l.$$
By using Propositions \ref{8.1} and \ref{8.2} we easily see that
$$R^{\Lambda}=U(\mathfrak{s}_{-})\spn\left\{R_{p}^{\Lambda}(n)v_{\Lambda} : p\geq l+1,\, n\leq -p\right\}.$$
\begin{center}\vspace{10pt}
\begin{tikzpicture}[scale=1]
\draw[->,thick] (-0.7,0) -- (5.5,0);
\draw[->,thick] (0,-0.7) -- (0,5.5);
\draw[dashed] (-0.1,2) -- (4.9,2);
\draw[dashed] (2,-0.1) -- (2,4.9);
\draw[dashed] (-0.1,-0.1) -- (4.9,4.9);
\filldraw (2.2,2.2) circle (3pt);\filldraw (2.2,2.7) circle (3pt);\filldraw (2.2,3.2) circle (3pt);\filldraw (2.2,3.7) circle (3pt);\filldraw (2.2,4.2) circle (3pt);\filldraw (2.2,4.7) circle (3pt);
\filldraw (2.7,2.7) circle (3pt);\filldraw (2.7,3.2) circle (3pt);\filldraw (2.7,3.7) circle (3pt);\filldraw (2.7,4.2) circle (3pt);\filldraw (2.7,4.7) circle (3pt);
\filldraw (3.2,3.2) circle (3pt);\filldraw (3.2,3.7) circle (3pt);\filldraw (3.2,4.2) circle (3pt);\filldraw (3.2,4.7) circle (3pt);
\filldraw (3.7,3.7) circle (3pt);\filldraw (3.7,4.2) circle (3pt);\filldraw (3.7,4.7) circle (3pt);
\filldraw (4.2,4.2) circle (3pt);\filldraw (4.2,4.7) circle (3pt);
\filldraw (4.7,4.7) circle (3pt);
\draw (5.5,-0.5) node {charge $p$};\draw (2,-0.5) node {$k/2$};
\draw (-1.2,5.5) node {degree $n$};\draw (-0.9,2) node {$-k/2$};
\end{tikzpicture}
\end{center}\vspace{10pt}
Since $p+q=k$ and $p>k/2$ implies $q<k/2$, we have for $n\leq -p$
$$R_{p}(n)v_{\Lambda}=a_pX^{(p)}(n)v_{\Lambda}\quad\text{mod }M_{(p-1)}$$
and
$$\lt R_{p}(n)v_{\Lambda}=X^{(p)}(n)v_{\Lambda}.$$
Note that in this case there is no $I_{p}^{\Lambda}(n)v_{\Lambda}$, i.e., initial conditions are ``trivial''.
Also note that 
$$\dim W(\Lambda)_{-l-1}=1\quad\textrm{and}\quad\dim W(\Lambda)_{-l-2}=3$$
and that $V''$ contains basis vectors
\begin{align*}
&X^{(l+1)}(-l-1)v_{\Lambda}\in V_{-l-1}''\quad\textrm{and}\\
&B(-1)X^{(l+1)}(-l-1)v_{\Lambda},\, X^{(l+1)}(-l-2)v_{\Lambda},\, X^{(l+2)}(-l-2)v_{\Lambda}\in V_{-l-2}''.
\end{align*} 
Hence $V''$ contains singular vectors
$f_{0}^{l +1}v_{\Lambda}$ and $f_{1}^{l +2}v_{\Lambda}$, and this implies
$$V''=W(\Lambda).$$
In order to reduce a spanning set of $W(\Lambda)=U(\mathfrak{n}_{-})R^{\Lambda}$
to a basis we need additional relations among relations.

\begin{pro}
Let $p\in \mathbb{Z}$. Then
\begin{align}
&\hspace{-174pt}\textrm{(a)}\,\, 2(p+1)R_{p+1}(\chi)=\lim_{\zeta\to\chi} (\chi/\zeta) (1+\zeta/\chi)^{2} R_{p}(\chi)X(\zeta),\label{XR:lim:1}\\
&\hspace{-174pt}\textrm{(b)}\,\, 8D_{\chi}R_{p+1}(\chi)=\lim_{\zeta\to\chi} D_{\zeta}(\chi/\zeta)^2 (1+\zeta/\chi)^{4} R_{p}(\chi)X(\zeta).\label{XR:lim:2}
\end{align}
\end{pro}

\begin{prf}
Note that (\ref{1:1}) implies that the limits in (\ref{XR:lim:1}) and (\ref{XR:lim:2}) exist, so we can define
\begin{align*}
&A_{p}(\chi)=\lim_{\zeta\to\chi} (\chi/\zeta) (1+\zeta/\chi)^{2} R_{p}(\chi)X(\zeta),\\
&B_{p}(\chi)=\lim_{\zeta\to\chi} D_{\zeta}(\chi/\zeta)^2 (1+\zeta/\chi)^{4} R_{p}(\chi)X(\zeta).
\end{align*}
Let 
$$A_p (\chi,\zeta)=\sum_{n,m\in\mathbb{Z}}A_{p}(n,m)\chi^{n}\zeta^{m}=(\chi/\zeta) (1+\zeta/\chi)^{2}R_p(\chi)X(\zeta).$$
Then for any $v\in M(\Lambda)$ and $N\in\mathbb{Z}$ the set 
$$\left\{(n,m)\in\mathbb{Z}^2: n+m=N,\,A_{p}(n,m)v\neq 0 \right\}$$
is finite. Set
$$f_{N}(\chi,\zeta)v=\sum_{n+m=N}A_{p}(n,m)v\chi^{n}\zeta^{m}\in M(\Lambda)[\chi^{\pm1},\zeta^{\pm 1}].$$
Note that
\begin{equation}\label{apovi}
A_{p}(\chi,\zeta)v=\sum_{N\in\mathbb{Z}}f_{N}(\chi,\zeta)v=\sum_{N\in\mathbb{Z}}\chi^{N}f_{N}(1,\zeta/\chi)v
\end{equation}
and
\begin{equation}\label{apovi2}
A_{p}(\chi)v=\lim_{\zeta\to\chi}A_{p}(\chi,\zeta)v=\sum_{N\in\mathbb{Z}}\chi^{N}f_{N}(1,1)v.
\end{equation}
Also note that
\begin{align}
B_{p}(\chi)v&=\lim_{\zeta\to\chi} D_{\zeta}(\chi/\zeta)(1+\zeta/\chi)^{2}A_{p}(\chi,\zeta)v\nonumber\\
&=\lim_{\zeta\to\chi}D_{\zeta}(\chi/\zeta)(1+\zeta/\chi)^2 \sum_{N\in\mathbb{Z}}\chi^{N}f_{N}(1,\zeta/\chi)v.\label{bpovi}
\end{align}
By using
$$\frac{z}{(1-z)^{2}}-\frac{z^{-1}}{(1-z^{-1})^{2}} = D\delta(z)$$
the left side of (\ref{2:2}) can be written as
\begin{equation}\label{fin:1}
D\delta(\zeta/\chi)(\chi/\zeta)\left(1+\zeta/\chi\right)^{2}R_{p}(\chi)X(\zeta)= D\delta(\zeta/\chi)A_p (\chi,\zeta).
\end{equation}
By applying operator $D$ on
$$\delta(z)f(z)=\delta(z)f(1),\quad f(z)\in\mathbb{C}[z^{\pm 1}],$$
we get
\begin{equation}\label{Ddelta}
D\delta(z)f(z)=D\delta(z)f(1)-\delta(z)Df(1),\quad f(z)\in\mathbb{C}[z^{\pm 1}].
\end{equation}
By using (\ref{apovi})--(\ref{Ddelta}) we get
\begin{align}
&(\chi/\zeta)(1+\zeta/\chi)^2 D\delta(\zeta/\chi)A_p (\chi,\zeta)v=(\chi/\zeta)(1+\zeta/\chi)^2\left(\sum_{r\in\mathbb{Z}}r(\zeta/\chi)^r\right) \left(\sum_{N\in\mathbb{Z}}f_{N}(\chi,\zeta)v\right)\nonumber\\
&=(\chi/\zeta)(1+\zeta/\chi)^2\sum_{N\in\mathbb{Z}} D\delta(\zeta/\chi) f_{N}(\chi,\zeta)v\nonumber\\
&=4D\delta(\zeta/\chi)\sum_{N\in\mathbb{Z}}  \chi^{N}f_{N}(1,1)v - \delta(\zeta/\chi)\lim_{\zeta\to\chi}D_{\zeta}(\chi/\zeta)(1+\zeta/\chi)^2 \sum_{N\in\mathbb{Z}}\chi^{N}f_{N}(1,\zeta/\chi)v\nonumber\\
&=4D\delta(\zeta/\chi)A_{p}(\chi)v - \delta(\zeta/\chi)B_{p}(\chi)v.\label{fin:2}
\end{align}
Finally, (\ref{2:2}), (\ref{fin:1}) and (\ref{fin:2}) imply
\begin{align}
&4D\delta(\zeta/\chi)A_{p}(\chi)v - \delta(\zeta/\chi)B_{p}(\chi)v\label{fin:3}\\
&=(\chi/\zeta)(1+\zeta/\chi)^2\left(2(p+1)D\delta(\zeta/\chi)R_{p+1}(\chi)v-2\delta(\zeta/\chi)D_{\chi}R_{p+1}(\chi)v\right)\label{fin:4}
\end{align}
By equating coefficients of $\zeta^0$ in (\ref{fin:3}) and (\ref{fin:4}) we get
$$-B_{p}(\chi)v=-8 D_{\chi}R_{p+1}(\chi)v,$$
and by equating coefficients of $\zeta^1$ in (\ref{fin:3}) and (\ref{fin:4}) we get
$$4\chi^{-1}A_{p}(\chi)v-\chi^{-1}B_{p}(\chi)v = -8\chi^{-1}D_{\chi}R_{p+1}(\chi)v + 8(p+1)\chi^{-1}R_{p+1}(\chi)v,$$
thus proving the proposition.
\end{prf}

The next two lemmas present  generalization of Lemma \ref{lemma2}, by providing us with relations among $R_{p}(\zeta)$ and $X^{(s)}(\xi)$.

\begin{lem}\label{lemma2_vol2}
Let $P(\zeta_1,\ldots,\zeta_s;\xi)$ be a Laurent polynomial, $s\leq p$, such that there exists a limit
$$\lim_{\zeta_i,\xi\to\chi} P(\zeta_1,\ldots,\zeta_s;\xi) \left(\prod_{i=1}^{s}D_{\zeta_i}^{n_i}X(\zeta_i)\right) R_{p}(\xi)$$
for $n_1+\ldots+n_s\leq 2s-1$, and that for some Laurent polynomial $R(\zeta,\xi)$ we have
$$\lim_{\zeta_i\to\zeta}P(\zeta_1,\ldots,\zeta_s;\xi) \left(\prod_{i=1}^{s}X(\zeta_i)\right) R_{p}(\xi)
=R(\zeta,\xi)X^{(s)}(\zeta)R_{p}(\xi).$$
Then for $r=0,1,\ldots,2s-1$
\begin{align}
\lim_{\zeta,\xi\to\chi}D_{\zeta}^{r}\left(R(\zeta,\xi)X^{(s)}(\zeta)R_{p}(\xi)\right)&\label{lemma2:1_vol2}\\
\equiv \left(R_{p+1}(\chi)\right),&\label{lemma2:2_vol2}
\end{align}
where (\ref{lemma2:2_vol2}) means that the coefficients of the formal Laurent series (\ref{lemma2:1_vol2}) in indeterminate $\chi$ are (``infinite'')
linear combinations of products of $X(n)'s$ with at least one coefficient of $R_{p+1}(\chi)$.
\end{lem}

\begin{prf}
Note that we can proceed in the same way as in the proof of Lemma \ref{lemma2}. Naturally, instead of using  (\ref{Xp})
in order to obtain (\ref{lemma2:2}), we will use (\ref{XR:lim:1}) and get (\ref{lemma2:2_vol2}).
\end{prf}

\begin{lem}\label{lemma2_vol3}
Let $P(\zeta_1,\ldots,\zeta_p;\xi_1,\ldots,\xi_s)$ be a Laurent polynomial, $s> p$, such that there exists a limit
$$\lim_{\zeta_i,\xi_j\to\chi} P(\zeta_1,\ldots,\zeta_p;\xi_1,\ldots,\xi_s) \left(\prod_{j=1}^{s}X(\xi_j)\right) \left(\prod_{i=1}^{p-1}D_{\zeta_i}^{n_i}X(\zeta_i)\right) D_{\zeta_p}^{n_p}R_{1}(\zeta_p)$$
for $n_1+\ldots+n_p\leq 2p-1$, and that for some Laurent polynomial $R(\zeta,\xi)$ we have
$$\lim_{\substack{\zeta_i\to\zeta\\\xi_j\to\xi}}P(\zeta_1,\ldots,\zeta_p;\xi_1,\ldots,\xi_s) \left(\prod_{j=1}^{s}X(\xi_j)\right)   \left(\prod_{i=1}^{p-1}X(\zeta_i)\right) R_{1}(\zeta_p)
=R(\zeta,\xi) X^{(s)}(\xi) R_p (\zeta).$$
Then for $r=0,1,\ldots,2p-1$
\begin{align}
\lim_{\zeta,\xi\to\chi}D_{\zeta}^{r}\left(R(\zeta,\xi)R_p (\zeta)X^{(s)}(\xi)\right)&\label{lemma2:1_vol3}\\
\equiv \left(X^{(s+1)}(\chi)\right) + \left(R_{s+1}(\chi)\right),&\label{lemma2:2_vol3}
\end{align}
where (\ref{lemma2:2_vol3}) means that the coefficients of the formal Laurent series (\ref{lemma2:1_vol3}) in indeterminate $\chi$ are (``infinite'')
linear combinations of products of $X(n)'s$ and (perhaps one) $R_{1}(m)$ with at least one coefficient of $X^{(s+1)}(\chi)$ or $R_{s+1}(\chi)$.
\end{lem}

\begin{prf}
First we derive a relation among
$R_{p'}(\zeta)$ and $X(\xi)$. 
Notice that
\begin{align}
&D_{\zeta} (\zeta/\xi)^2 (1+\xi/\zeta)^{4}=2 (\zeta/\xi) (1+\xi/\zeta)^{2} (-\xi/\zeta + \zeta/\xi).\label{8:8:temp}\\
&D_{\xi} (\zeta/\xi)^2 (1+\xi/\zeta)^{4}=2 (\zeta/\xi) (1+\xi/\zeta)^{2} (\xi/\zeta - \zeta/\xi).\label{8:9:temp}
\end{align}
By applying $D_{\chi}$ on 
$$8(p'+1)R_{p'+1}(\chi)=\lim_{\zeta,\xi\to\chi} (\zeta/\xi)^2 (1+\xi/\zeta)^{4} R_{p'}(\zeta)X(\xi)$$
(recall (\ref{XR:lim:1})) and using (\ref{XR:lim:2}) and (\ref{8:8:temp}) we get
\begin{align*}
&8(p'+1)D_{\chi}R_{p'+1}(\chi)=D_{\chi}\lim_{\zeta,\xi\to\chi} (\zeta/\xi)^2 (1+\xi/\zeta)^{4} R_{p'}(\zeta)X(\xi)\\
&=\lim_{\zeta,\xi\to\chi} (D_{\zeta} + D_{\xi})(\zeta/\xi)^2 (1+\xi/\zeta)^{4} R_{p'}(\zeta)X(\xi)\\
&=\lim_{\zeta,\xi\to\chi} \bigg( 2 (\zeta/\xi) (1+\xi/\zeta)^{2} (-\xi/\zeta + \zeta/\xi)R_{p'}(\zeta)X(\xi)
 + (\zeta/\xi)^2 (1+\xi/\zeta)^{4}D_{\zeta}R_{p'}(\zeta)X(\xi)\\
&\hspace{48pt}+ D_{\xi}(\zeta/\xi)^2 (1+\xi/\zeta)^{4}R_{p'}(\zeta)X(\xi)\bigg)\\
 &=\lim_{\zeta,\xi\to\chi} \bigg( (\zeta/\xi)^2 (1+\xi/\zeta)^{4}D_{\zeta}R_{p'}(\zeta)X(\xi)+D_{\xi}(\zeta/\xi)^2 (1+\xi/\zeta)^{4}R_{p'}(\zeta)X(\xi)\bigg)\\
 &=\lim_{\zeta,\xi\to\chi} (\zeta/\xi)^2 (1+\xi/\zeta)^{4}D_{\zeta}R_{p'}(\zeta)X(\xi) + 8D_{\chi}R_{p'+1}(\chi),
\end{align*}
i.e.
\begin{equation}\label{8:10:finally}
8p'D_{\chi}R_{p'+1}(\chi)=\lim_{\zeta,\xi\to\chi} (\zeta/\xi)^2 (1+\xi/\zeta)^{4}D_{\zeta}R_{p'}(\zeta)X(\xi).
\end{equation}
Finally, (\ref{XR:lim:2}), (\ref{8:9:temp}) and (\ref{8:10:finally}) imply
\begin{equation}\label{a_relation}
\lim_{\zeta,\xi\to\chi} (\zeta/\xi)^2 (1+\xi/\zeta)^{4}D_{\zeta}R_{p'}(\zeta)X(\xi)=p'\lim_{\zeta,\xi\to\chi} (\zeta/\xi)^2 (1+\xi/\zeta)^{4}R_{p'}(\zeta)D_{\xi}X(\xi).
\end{equation}

Now we can proceed in the same way as in the proof of Lemma \ref{lemma2}: relation (\ref{a_relation}) allows us to obtain (\ref{lemma2:2_vol3}).
More precisely, 
for $0\leq r\leq 2p-1$ we have
 \begin{equation*}
 (D_{\zeta_1}+\ldots+D_{\zeta_p})^{r}=\sum_{i_1,\ldots,i_r}D_{\zeta_{i_1}}\cdots D_{\zeta_{i_r}},
 \end{equation*}
 where for each choice of $(i_1,\ldots,i_r)$ there is an index $t\in\left\{1,\ldots,p\right\}$ such that there is at most one $D_{\zeta_t}$  in the product
 $D_{\zeta_{i_1}}\cdots D_{\zeta_{i_r}}$, 
 so the formal Laurent series
 $$(D_{\zeta_1}+\ldots+D_{\zeta_p})^{r}P(\zeta_1,\ldots,\zeta_p;\xi_1,\ldots,\xi_s) \left(\prod_{j=1}^{s}X(\xi_j)\right)  \left(\prod_{i=1}^{p-1}X(\zeta_i)\right) R_{1}(\zeta_p)$$
 can be written as a sum over $i_1,\ldots,i_r$. Any such summand can be written as
  $$\left(\prod_{i\neq t}D_{\zeta_i}^{a_i}\right)D_{\zeta_t}^{a}P(\zeta_1,\ldots,\zeta_p;\xi_1,\ldots,\xi_s) \left(\prod_{j=1}^{s}X(\xi_j)\right)  \left(\prod_{i=1}^{p-1}X(\zeta_i)\right)R_{1}(\zeta_p),$$
  where $a=0,1$. 
  If $t<p$ we obtain, as in the proof of Lemma \ref{lemma2}, a summand containing a coefficient of $X^{(s+1)}(\chi)$. 
  If $t=p$ and $a=0$ we obtain, by using (\ref{XR:lim:1}), a summand containing a coefficient of $R_{s+1}(\chi)$.
  If $t=p$ and $a=1$ we obtain, by using (\ref{a_relation}), (\ref{XR:lim:1}) and (\ref{XR:lim:2}), a summand containing a coefficient of $R_{s+1}(\chi)$.
 \end{prf}
 
 Next, we have a generalization of Lemma \ref{lemma5}:
 
 \begin{lem}\label{lemma5_vol2}
Let $r=\min\left\{s,p\right\}$, $q=\max\left\{s,p\right\}$  and $n,j\in\mathbb{Z}$. For any given vector $v$ in any given highest weight $\hat{\mathfrak{g}}$-module $V$ 
denote by $R^{(s,p)}_{j,n}v$ a sequence of $2r$ vectors
\begin{align}
&X^{(s)}(j)R_{p}(n-j)v,\nonumber\\
&X^{(s)}(j+1)R_{p}(n-(j+1))v,\label{lemma5:summands_vol2}\\
&\qquad\qquad\qquad\vdots\nonumber\\
&X^{(s)}(j+2r-1)R_{p}(n-(j+2r-1))v\nonumber
\end{align}
in a formal Laurent series $X^{(s)}(\zeta)R_{p}(\xi)v$. Then each vector in $R^{(s,p)}_{j,n}v$
can be written as a linear combination of vectors
$$X^{(s)}(i)R_{p}(n-i)v\notin R^{(s,p)}_{j,n}v$$
and vectors obtained by the action of monomials, that have  a coefficient of $X^{(q+1)}(\chi)$ or $R_{q+1}(\chi)$ as a factor, on $v$.
\end{lem}

\begin{prf}
By using Lemmas \ref{lemma2_vol2} for $s\leq p$ and \ref{lemma2_vol3} for $s>p$ we can prove the lemma in a way analogous to the proof 
of Lemma \ref{lemma5}. First, we notice that an analogue of Lemma \ref{lemma4} holds in both cases.
Next, we see that, when we try to express vectors (\ref{lemma5:summands_vol2})
as a certain linear combination described above, we will, as before, get coefficients $a_{r,m}$ that form a regular matrix (recall Lemma \ref{lemma3}).
\end{prf}

The main result in this section is

\begin{thm}\label{last:thm}
Let $\Lambda=(l+1)\Lambda_0 + l\Lambda_1$ as in (\ref{lambda}). The set $\mathcal{B}_{W(\Lambda)}$ is a basis
of the maximal submodule $W(\Lambda)$ of Verma module $M(\Lambda)$.
\end{thm}

\begin{prf}
By using lemmas \ref{lemma5_vol2} and \ref{lemma5} we can prove that the set of vectors of the form (\ref{setofvectors}) spans $V''$.
Furthermore, for $D'=D_{p'}^{\Lambda}(n')$, $D''=D_{p''}^{\Lambda}(n'')$ and  
$$x=B(i_1)\cdots B(i_r)X^{(p_1)}(j_1)\cdots X^{(p_s)}(j_s)$$
we get, by using 
Proposition \ref{pro:Xp:B},
\begin{align*}
& xD'D''v_{\Lambda}=x(\lt(D')\textrm{ }+\textrm{ higher terms})D''v_{\Lambda}=x\lt(D')D''v_{\Lambda}\textrm{ }+\textrm{ higher terms},\\
& xD'D''v_{\Lambda}=xD'(\lt(D'')\textrm{ }+\textrm{ higher terms})v_{\Lambda}=xD'\lt(D'')v_{\Lambda}\textrm{ }+\textrm{ higher terms},
\end{align*}
so we conclude that two vectors of the form (\ref{setofvectors}) with a same leading term are proportional modulo higher terms in our order.
\end{prf}

As a consequence we have a vertex-operator theoretic proof of Andrews' analytic  identity
\begin{equation*}
\prod_{\substack{n\geq 1\\n\not\equiv 0,\pm (l+2)(\textrm{mod }2l+3)}}(1-q^n)^{-1} = 
\sum_{n_1 , n_2 , \ldots ,  n_{l}\geq 0} \frac{ q^{ N_{1}^{2}+N_{2}^{2}+...+ N^{2}_{l} } }     { (q)_{n_1 } (q)_{n_2 }\cdots (q)_{n_{l}}}.
\end{equation*}

\section*{Acknowledgement}
The authors would like to thank Jim Lepowsky for his  valuable comments and suggestions.


\end{document}